\documentclass[12pt]{article}

\usepackage{a4wide}
\usepackage{amssymb}
\usepackage{amsfonts}
\usepackage{amsmath}
\input xy
\xyoption{arrow} \xyoption{matrix}

\date{}

\newtheorem{proposition}{Proposition}[section]
\newtheorem{theorem}[proposition]{Theorem}
\newtheorem{lemma}[proposition]{Lemma}

\newtheorem{corollary}[proposition]{Corollary}

\def\GK{{\rm  GK}\,}

\def\Hom{{\rm Hom}}
\def\der{\partial }

\def\nFM0{{\nu }_{F,M_0}}
\def\nFN0{{\nu }_{F,N_0}}
\def\nGN0{{\nu }_{G,N_0}}

\def\N0{ {\bf N}_0 }

\def\t{\otimes}
\def\g{\gamma}
\def\v{\varphi}
\def\ra{\rightarrow}

\def\Xpm{X^{\pm }}

\def\s{\sigma}
\def\Z{\mathbb{Z}}

\def\l1{{\lambda}_1}

\def\a{\alpha}
\def\a0{ {\alpha }_0}
\def\a1{ {\alpha }_1}

\def\l{\lambda}

\def\nFGM0{{\nu }_{F,G,M_0}}

\def\nFN0{{\nu}_{F,N_0}}

\def\sm{{\sigma}^m}

\def\sm1{{\sigma}^{-1}}

\def\smtp1{{\sigma}^{-t+1}}

\def\S1{S^{-1}}

\def\Xpm1{X^{\pm 1}_1}

\def\sPM1{{\sigma }^{\pm 1}}
\def\sMP1{{\sigma }^{\mp 1 }}

\def\d{\delta}

\def\di{{\rm d.ind}}

\def\L{\Lambda}

\def\CA{{\cal A}}

\def\CD{{\cal D}}

\def\Ytm1{Y^{t-1}}
\def\Yim1{Y^{i-1}}

\def\CM{{\cal M}}
\def\CN{{\cal N}}

\def\CF{{\cal F}}
\def\CG{{\cal G}}
\def\CH{{\cal H}}

\def\supp{{\rm supp }}

\def\ker{ {\rm ker } }

\def\gr{ {\rm gr} }

\def\SL2Z{ {\rm SL}_2({\bf Z}) }

\def\th{ \theta }

\def\Gp1{ G^{1 , 1 } }
\def\P11{ P^{-1 , 1 } }
\def\Pp1{ P^{1 , 1 } }

\def\th{\theta}

\def\nCLsr{{}^\nu\kern-2pt {\cal L}^{\sigma , \rho  }}
\def\nP{{}^\nu \kern-2pt P}
\def\nL{{}^\nu\kern-2pt L}
\def\nLL{{}^\nu\kern-2pt \Lambda}
\def\nPsr{{}^\nu\kern-2pt P^{\sigma , \rho  }}
\def\nLsr{{}^\nu\kern-2pt L^{\sigma , \rho  }}
\def\nuCL{{}^\nu\kern-2pt  {\cal L}}
\def\nCLsr{{}^\nu\kern-2pt {\cal L}^{\sigma , \rho  }}
\def\nCL1m{{}^\nu\kern-2pt {\cal L}^{-1 , 1  }}
\def\x1nu{x^\frac{1}{\nu}}
\def\xm1nu{x^{-\frac{1}{\nu}}}

\def\CN{{\cal N}}
\def\ra{\rightarrow }

\def\CB{{\cal B}}

\def\CI{{\cal I}}

\def\CC{ {\cal C}}

\def\CH{ {\cal H}}
\def\CP{ {\cal P}}

\def\nAM0{{\nu }_{{\cal A},M_0}}
\def\nAN0{{\nu }_{{\cal A},N_0}}

\def\End{ {\rm End }}

\def\CP{ {\cal P }}
\def\det{ {\rm det }}

\def\ga{\mathfrak{a}}
\def\gb{\mathfrak{b}}
\def\gc{\mathfrak{c}}

\def\gm{\mathfrak{m}}
\def\gp{\mathfrak{p}}
\def\gq{\mathfrak{q}}
\def\gr{\mathfrak{r}}

\def\GL{{\rm GL}}
\def\SL{{\rm SL}}

\def\Spec{{\rm Spec}}

\def\Hom{{\rm Hom}}

\def\di!{\frac{\der^i}{i!}}
\def\dik!{\frac{\der^k_i}{k!}}

\def\Max{{\rm Max}}

\def\N{\mathbb{N}}

\def\0{\overline{0}}
\def\1{\overline{1}}

\def\Ln1{\L_{n,\overline{1}}}

\def\a1{a_{\overline{1}}}

\def\S{\Sigma}

\def\vn1{\overrightarrow{n-1}}

\def\hx{\widehat{x}}

\def\Q{\mathbb{Q}}
\def\im{{\rm im}}

\def\mA{\mathbb{A}}
\def\mD{\mathbb{D}}
\def\mL{\mathbb{L}}

\def\gf{\mathfrak{f}}
\def\soc{{\rm soc}}
\def\Sub{{\rm Sub}}
\def\SSub{{\rm SSub}}
\def\Inc{{\rm Inc}}
\def\Min{{\rm Min}}
\def\csupp{{\rm csupp}}

\begin{document}

\author{V. V. \  Bavula 
}

\title{The  Jacobian  algebras }

\maketitle
\begin{abstract}
Let $P_n:=K[x_1, \ldots , x_n]$ be a polynomial algebra over a
field $K$ of characteristic zero. The {\em Jacobian algebra}
$\mA_n$ is the subalgebra of $\End_K(P_n)$ generated by the Weyl
algebra $A_n:= \CD (P_n)=K\langle x_1, \ldots , x_n, \der_1,
\ldots , \der_n\rangle$ and the elements $(\der_1x_1)^{-1}, \ldots
, (\der_nx_n)^{-1}\in \End_K(P_n)$. The algebra $\mA_n$ appears
naturally in study of the group of  automorphisms of $P_n$. The
 algebra $\mA_n$ is large since  it contains a ring
 $M_\infty (K):=\varinjlim M_d(K)$ of infinite dimensional matrices
 and all (formal)
integro-differential operators as all
$\int_i=x_i(\der_ix_i)^{-1}\in \mA_n$. Surprisingly, the algebras
$\mA_n$ and $A_n$ have little in common: the algebra $\mA_n$ is
neither left nor right Noetherian (even contains infinite direct
sums of nonzero left and right ideals); not simple; not a domain;
contains nilpotent elements; local;  prime; central; self-dual;
$\GK (\mA_n) = 3n$; ${\rm cl.K.dim} (\mA_n) = n$; has only
finitely many, say $s_n$,  ideals ($2-n+\sum_{i=1}^n2^{n\choose
i}\leq s_n \leq 2^{2^n}$) which are found explicitly (they
commute, $IJ=JI$; each of them is an  idempotent ideal, $I^2=I$,
and $\GK (\mA_n/I)=3n$ if $I\neq \mA_n$). $\Spec (\mA_n )$ is
found, it contains exactly $2^n$ elements, each nonzero prime
ideal is a unique sum of primes of height 1 (and any such a sum is
a prime ideal). Each nonzero ideal is a {\em unique} product and a
{\em unique} intersection of incomparable primes; moreover in both
presentations the primes are the same, they are the {\em minimal}
primes over the ideal.  The group of units $\mA_n^*$ of $\mA_n$ is
huge ($\mA_n^*\supseteq K^*\times ( (\Z^n)^{(\Z )}\ltimes
\GL_\infty (K))$).  $\mA_n$ has only one {\em faithful simple}
module, namely $P_n$.



$${\bf Contents}$$
\begin{enumerate}
\item Introduction. \item The Jacobian algebras and localizations
of the Weyl algebras. \item Unique factorization of ideals of
$\mA_n$ and $\Spec (\mA_n)$. \item The group of units $\mA_n^*$ of
$\mA_n$.
\end{enumerate}
\end{abstract}


\section{Introduction}
Throughout, ring means an associative ring with $1$. Let $K$ be a
commutative $\Q$-algebra, $K^*$ be its group of units,  $P_n:=
K[x_1, \ldots , x_n]$ be a polynomial algebra over $K$;
$\der_1:=\frac{\der}{\der x_1}, \ldots , \der_n:=\frac{\der}{\der
x_n}$ be the partial derivatives ($K$-linear derivations) of
$P_n$.

{\bf General properties of the Jacobian algebras}. The {\em
Jacobian algebra} $\mA_n$ is the subalgebra of $\End_K(P_n)$
generated by the Weyl algebra $A_n:=K\langle x_1, \ldots , x_n,
\der_1, \ldots , \der_n\rangle$ and the elements $H_1^{-1}, \ldots
, H_n^{-1}\in \End_K(P_n)$ where $H_1:= \der_1x_1, \ldots , H_n:=
\der_nx_n$. Clearly, $\mA_n = \mA_1(1) \t \cdots \t \mA_1(n)
\simeq \mA_1^{\t n }$ where $\mA_1(i) := K\langle x_i, \der_i ,
H_i^{-1}\rangle \simeq \mA_1$. The algebra $\mA_n$ contains all
the  integrations $\int_i: P_n\ra P_n$, $ p\mapsto \int p \,
dx_i$, since  $\int_i= x_iH_i^{-1}$. In particular, the algebra
$\mA_n$ contains all (formal) integro-differential operators with
polynomial coefficients. This fact explains $(i)$ significance of
the algebras $\mA_n$ for Algebraic Geometry and the theory of
integro-differential operators; $(ii)$ why the algebras $\mA_n$
and $A_n$ have different properties, and $(iii)$ why the group
$\mA_n^*$ of units of the algebra $\mA_n$ is huge (there are many
invertible integro-differential operators).

Till the end of this section,  $K$ is a field of characteristic
zero. When $n=1$ the group $\mA_1^*$ is found explicitly,
$\mA_1^*\simeq K^*\times ( \Z^{(\Z )}\ltimes \GL_\infty (K))$
 (Theorem \ref{18Ma7}) as well as an inversion formula $u^{-1}$
 for $u\in \mA_1^*$. This gives explicitly polynomial solutions
 for {\em all invertible}  integro-differential operators on an affine
 line: $uy= f$ $\Rightarrow$ $y= u^{-1} f$ where $f\in K[x_1]$ and
 $y$ is an unknown. For $n\geq 2$, a description of the group
 $\mA_n^*$ is given (Theorem \ref{24Ma7}), it looks like it is a
 challenging problem to find an inversion formula for $u\in \mA_n^*$ (one should go
 far beyond the Dieudonn\'{e} determinant). Though, a criterion of
 invertibility is found (Theorem \ref{23Ma7}). Moreover, the group
 $\mA_n^*$ contains the subgroup $ K^*\times ( (\Z^n)^{(\Z )}\ltimes
\GL_\infty (K))$ elements of which  are called {\em minimal
integro-differential operators}. For each such an operator $u$ one
can write down an inversion formula $u^{-1}$ in the same manner as
in the case $n=1$, and, therefore, one obtains explicitly
polynomial solutions for all minimal integro-differential
equations $uy=f$ where $f\in P_n$.

The Weyl algebra $A_n= A_n(K)$ is a simple, Noetherian domain of
Gelfand-Kirillov dimension $\GK (A_n) =2n$. The Jacobian algebra
$\mA_n$ is neither  left nor right Noetherian, it  contains
infinite direct sums of nonzero left and right ideals. This means
that the concept of the left (and right) Krull dimension makes no
sense for $\mA_n$ but the classical Krull dimension of $\mA_n$ is
$n$ (Corollary \ref{k28Ma7}). The algebra $\mA_n$ is a central,
prime algebra of Gelfand-Kirillov dimension $3n$ (Corollary
\ref{15Ma7}).

The canonical involution $\th$ of the Weyl algebra can be extended
to the algebra $\mA_n$ (see (\ref{thinv})). This means that the
algebra $\mA_n$ is {\em self-dual} ($\mA_n\simeq \mA_n^{op}$), and
so its left and right algebraic  properties are the same.  Note
that the {\em Fourier transform} on the Weyl algebra $A_n$ can not
be lifted to $\mA_n$. Many properties of the algebra
$\mA_n=\mA_1^{\t n }$ are determined by properties of $\mA_1$.
When $n=1$ we usually drop the subscript `1' in $x_1$, $\der_1$,
$H_1$, etc. The algebra $\mA_1$ contains the only proper ideal
$F=\oplus_{i,j\in \N} KE_{ij}$ where
$$
E_{ij}:=\begin{cases}
x^{i-j}(x^j\frac{1}{\der^jx^j}\der^j-x^{j+1}\frac{1}{\der^{j+1}x^{j+1}}\der^{j+1})& \text{if $i\geq j$},\\
(\frac{1}{\der x}\der )^{j-i}(x^j\frac{1}{\der^jx^j}\der^j-x^{j+1}\frac{1}{\der^{j+1}x^{j+1}}\der^{j+1})& \text{if $i<j$}.\\
\end{cases}
$$
As a ring without 1, the ring $F$ is canonically isomorphic to the
ring $M_\infty (K):=\varinjlim M_d(K)= \oplus_{i,j\in \N} KE_{ij}$
of infinite-dimensional matrices where $E_{ij}$ are the matrix
units ($F\ra M_\infty  (K)$, $E_{ij}\mapsto E_{ij}$). This is a
very important fact as we can apply  concepts of
finite-dimensional linear algebra (like trace, determinant, etc)
to integro-differential operators which is not obvious from the
outset. This fact is crucial in finding an inversion formula for
 elements of $\mA_1^*$.

The algebra $\mA_n = \oplus_{\alpha \in \Z^n} \mA_{n, \alpha}$ is
a $\Z^n$-graded algebra where $\mA_{n,\alpha}:=\t_{k=1}^n
\mA_{1,\alpha_k}(k)$ and, for $n=1$, (Theorem \ref{8Ap7})
$$
\mA_{1,i}=\begin{cases}
x^i\mD_1& \text{if $i\geq 1$},\\
\mD_1& \text{if $i=0$},\\
\mD_1 \der^{-i}& \text{if $i\leq -1$},\\
\end{cases}
$$
where $\mD_1:= L\oplus (\oplus_{i,j\geq 1}Kx^iH^{-j}\der^i)$ is a
{\em commutative, non-Noetherian} algebra and $L= K[H^{\pm 1},
(H+1)^{-1}, (H+2)^{-1}, \ldots ]$. This gives a `compact'
$K$-basis for the algebra $\mA_1$ (and $\mA_n$). This basis
`behaves badly' under multiplication. A more conceptual
(`multiplicatively friendly') basis is given in Theorem
\ref{12Ma7}.

\begin{itemize}
\item (Corollary \ref{15Ma7}.(10)) $P_n$ {\em is the only
faithful, simple $\mA_n$-module.}
\end{itemize}

{\bf Spec ($\mA_n$)}. $0$ is a prime ideal of $\mA_n$.
$$ \gp_1:=F\t\mA_{n-1}, \gp_2:= \mA_1\t F\t \mA_{n-2} , \ldots ,
\gp_n :=\mA_{n-1}\t F,$$  are precisely the prime ideals  of
height 1 of $\mA_n$. Let $\Sub_n$ be the set of all subsets of $\{
1, \ldots , n\}$.
\begin{itemize}
\item (Corollary \ref{c27Ma7}) {\em The map $\Sub_n\ra \Spec
(\mA_n)$, $ I\mapsto \gp_I:= \sum_{i\in I}\gp_i$, $\emptyset
\mapsto 0$, is a bijection, i.e. any nonzero prime ideal of
$\mA_n$ is a unique sum of primes of height 1; $|\Spec
(\mA_n)|=2^n$; the height of $\gp_I$ is $| I|$;  and}
 \item (Lemma \ref{p28Ma7}) $\gp_I\subseteq \gp_J$ {\em iff} $I\subseteq
 J$.
\item (Corollary \ref{b28Ma7}) $\ga_n:= \gp_1+\cdots +\gp_n$ {\em
is the only prime ideal of $\mA_n$ which is  completely prime;
$\ga_n$ is the only ideal $\ga$ of $\mA_n$ such that $\ga \neq
\mA_n$ and $\mA_n/\ga$ is a Noetherian (resp. left Noetherian,
resp. right Noetherian) ring.}
\end{itemize}

{\bf Ideals of $\mA_n$ and their unique factorization}. The ideal
theory of $\mA_n$ is `very arithmetic.' Let $\CB_n$ be the set of
all functions $f:\{ 1, 2, \ldots , n\} \ra  \{ 0,1\}$. For each
function $f\in \CB_n$, $I_f:= I_{f(1)}\t \cdots \t I_{f(n)}$ is
the ideal of $\mA_n$ where $I_0:=F$ and $I_1:= \mA_1$.  Let
$\CC_n$ be the set of all  subsets of $\CB_n$ all distinct
elements of which are incomparable (two distinct elements $f$ and
$g$ of $\CB_n$ are {\em incomparable} if neither $f(i)\leq  g(i)$
nor $f(i)\geq g(i)$ for all $i$). For each $C\in \CC_n$, let
$I_C:= \sum_{f\in C}I_f$, the ideal of $\mA_n$.  The next result
classifies all the ideals of $\mA_n$.

\begin{itemize}
\item (Theorem \ref{C15Ma7}) {\em The map $C\mapsto I_C:=
\sum_{f\in C}I_f$ from the set $\CC_n$ to the set of  ideals of
$\mA_n$ is a bijection where $I_\emptyset :=0$. In particular,
there are only finitely many ideals, say $s_n$,  of $\mA_n$.
Moreover,} $2-n+\sum_{i=1}^n2^{n\choose i}\leq s_n \leq 2^{2^n}$
(Corollary \ref{b27Ma7}). \item {\em Each ideal $I$ of $\mA_n$ is
an idempotent ideal, i.e. $I^2= I$. \item Ideals of $\mA_n$
commute} ($IJ= JI$). \item (Theorem \ref{28Ma7}) {\em The lattice
of ideals of $\mA_n$ is distributive}. \item (Corollary
\ref{15Ma7}.(4,7)) {\em The ideal $\ga_n$ is the largest (hence,
the only maximal) ideal of $\mA_n$ distinct from $\mA_n$,  and
$F^{\t n }$ is the smallest nonzero ideal of $\mA_n$. } \item
(Corollary \ref{15Ma7}.(11)) $\GK (\mA_n/ \ga )= 3n$ {\em for all
ideals $\ga$ of $\mA_n$ such that $\ga \neq \mA_n$. }
\end{itemize}

For each ideal $\ga$ of $\mA_n$, $\Min (\ga )$ denotes the set of
minimal primes over $\ga$. Two distinct prime ideals $\gp $ and
$\gq$ are called {\em incomparable} if neither $\gp \subseteq \gq$
nor $\gp\supseteq \gq$. The algebras $\mA_n$ have  beautiful ideal
theory  as the following unique factorization properties
demonstrate.

\begin{itemize}
\item (Theorem \ref{27Ma7}) 1. {\em Each ideal $\ga$ of $\mA_n$
such that $\ga \neq \mA_n$ is a unique product of incomparable
primes, i.e. if $\ga = \gq_1\cdots \gq_s= \gr_1\cdots \gr_t$ are
two such products then $s=t$ and $\gq_1= \gr_{\s (1)}, \ldots ,
\gq_s= \gr_{\s (s)}$ for a permutation $\s$ of $\{ 1, \ldots,
n\}$.}

2.  {\em Each ideal $\ga$ of $\mA_n$ such that $\ga \neq \mA_n$ is
a unique intersection of incomparable primes, i.e. if $\ga =
\gq_1\cap \cdots\cap  \gq_s= \gr_1\cap \cdots \cap \gr_t$ are two
such intersections then $s=t$ and $\gq_1= \gr_{\s (1)}, \ldots ,
\gq_s= \gr_{\s (s)}$ for a permutation $\s$ of $\{ 1, \ldots,
n\}$.}

3.  {\em For each ideal $\ga$ of $\mA_n$ such that $\ga \neq
\mA_n$, the sets of incomparable primes in statements 1 and 2 are
the same, and so $\ga=\gq_1\cdots \gq_s=\gq_1\cap \cdots\cap
\gq_s$.}

4. {\em  The ideals $\gq_1,\ldots , \gq_s$ in statement 3 are the
minimal primes of $\ga$, and so  $\ga = \prod_{\gp \in \Min (\ga
)}\gp =\cap_{\gp \in \Min (\ga )}\gp$.}
\end{itemize}

\begin{itemize}
\item (Corollary \ref{a28Ma7}) $\ga\cap \gb = \ga \gb$ {\em for
all ideals $\ga$ and $\gb$ of $\mA_n$.}
\end{itemize}
The next theorem gives all decompositions of an ideal as a product
or intersection of ideals.
\begin{itemize}
\item (Theorem  \ref{1Jun7}) {\em Let $\ga$ be an ideal of
$\mA_n$, and $\CM$ be the minimal elements with respect to
inclusion of a set of ideals $\ga_1, \ldots , \ga_k$ of $\mA_n$.
Then }

1. $\ga = \ga_1\cdots \ga_k$ {\em iff} $\Min (\ga ) = \CM$.

2. $\ga = \ga_1\cap \cdots\cap \ga_k$ {\em iff} $\Min (\ga ) =
\CM$.
\end{itemize}
This is a rare example of a non-commutative  algebra of Krull
dimension $>1$ where one has a complete picture of decompositions
of ideals.

{\bf The group $\mA_1^*$ of units  of $\mA_1$}. For each integer
$i\geq 1$, consider the element of $\mA_1^*$:
$$ (H-i)_1^{-1}:=
\begin{cases}
x\frac{1}{H^2}\der +1-x\frac{1}{H}\der& \text{if $i=1$ },\\
x\frac{1}{H\der^ix^i}\der+\sum_{j=0}^{i-2} \frac{1}{j+1-i}\, \pi_j +\pi_{i-1}& \text{if $i\geq 2$ },\\
\end{cases}
$$
where $\pi_j:=
x^j\frac{1}{\der^jx^j}\der^j-x^{j+1}\frac{1}{\der^{j+1}x^{j+1}}\der^{j+1}$.
Consider the following subgroup of $\mA_1^*$,
$$\CH := \{ \prod_{i\geq 0} (H+i)^{n_i}\cdot \prod_{i\geq
1}(H-i)^{n_{-i}}_1\, | \, (n_i)\in \Z^{(\Z )}\}\simeq \Z^{(\Z )}
$$
 where $\Z^{(\Z )}$ is the  direct sum of $\Z$
copies of the group $\Z$, see (\ref{defH1}) for detail. Let
$\GL_\infty (K):=\{ u\in 1+M_\infty (K)\, | \, \det (u)\neq 0\}$.
The group $(1+F)^*$ of units of the multiplicative monoid $1+F$ is
equal to $(1+F)^* = (1+M_\infty (K))^*= \GL_\infty (K)$. Note that
$(1+F)^*\subseteq \mA_1^*$.

\begin{itemize}
 \item (Theorem \ref{18Ma7}) 1. {\em $\mA_1^*= K^*\times (\CH \ltimes (1+F)^*)$, each unit $a$
of $\mA_1$ is a unique product $a= \l \alpha (1+f)$ for some
elements $\l \in K^*$, $\alpha \in \CH$, and $f\in F$ such that
$\det (1+f)\neq 0$.}

2.  $\mA_1^*= K^*\times (\CH \ltimes \GL_\infty (K))$.

3.  {\em The centre of the group $\mA_1^*$ is $K^*$.}

4. {\em The commutant $\mA_1^{*(2)}:= [\mA_1^*, \mA_1^*]$ of the
 group $\mA_1^*$ is equal to $\SL_\infty (K):=\{ v\in (1+F)^*=
 M_\infty (K)\, | \, \det (v) = 1\}$, and $\mA_1/[\mA_1^*,
 \mA_1^*]\simeq K^*\times \CH \times K^*$.}

5.  {\em All the higher commutants $\mA_1^{*(i)}:= [\mA_1^*,
 \mA_1^{*(i-1)}]$, $i\geq 2$, are equal to $\mA_1^{*(2)}$.}
\end{itemize}

{\bf The group of units $\mA_n^*$ of $\mA_n$}.

\begin{itemize}
 \item (Theorem \ref{24Ma7}) 1.  $\mA_n^*= K^*\times (\CH_n\ltimes (1+\ga_n)^*)$
 where $\CH_n:= \prod_{i=1}^n\CH (i)$ and $\ga_n:= \gp_1+\cdots +\gp_n$.

 2. {\em The
centre of the group $\mA_n^*$ is $K^*$.}
\end{itemize}

Theorem \ref{23Ma7} is a criterion of when an element of the
monoid $1+\ga_n$ belongs to its group  $(1+\ga_n )^*$ of units.

{\em Question. Is the global dimension of $\mA_n$ equal to $2n$
(or $\infty$)?}


\section{The Jacobian algebras and localizations
of the Weyl algebras}\label{WAIL}

In this section, two $K$-bases for the Jacobian algebras are found
(Theorems \ref{8Ap7} and \ref{12Ma7}), and several properties of
 the algebras $\mA_n$ are proved: $\mA_n$ is a central, prime,
self-dual,  non-Noetherian algebra.

We start by recalling some properties of generalized Weyl
algebras. Some of these algebras are factor algebras of the
Jacobian algebras.

{\bf Generalized Weyl Algebras}. Let $D$ be a ring,
$\sigma=(\sigma_1,...,\sigma_n)$ be
 an $n$-tuple  of  commuting automorphisms of $D$,  and
$a=(a_1,...,a_n)$ be  an $n$-tuple  of (non-zero) elements of the
centre $Z(D)$  of $D$ such that $\sigma_i(a_j)=a_j$ for all $i\neq
j$.

        The {\it generalized Weyl algebra} $A=D(\sigma,a)$ (briefly GWA) of
degree $n$ with the {\it base } ring  $D$  is  a  ring generated
by $D$  and    $2n$ indetermina\-tes $x_1,...,x_n,$ $y_1,...,$
$y_n$ subject to the defining relations \cite{Bav-FA91},
\cite{Bav-AlgAnaliz92}:
\begin{align*}
y_ix_i&=a_i,& x_iy_i&=\sigma_i(a_i), \\
x_i\alpha &=\sigma_i(\alpha)x_i, &
 y_i\alpha&=\sigma_i^{-1}(\alpha)y_i, \;\;\;  \alpha \in D,
\end{align*}
$$[x_i,x_j]=[y_i,y_j]=[x_i,y_j]=0, \;\;\;  i\neq j,$$
where $[x, y]=xy-yx$. We say that  $a$  and $\sigma $ are the sets
of {\it defining } elements and automorphisms of $A$ respectively.
The GWAs are also known as {\em hyperbolic rings}, see the book of
Rosenberg \cite{Rosen-book}.
 For a
vector $k=(k_1,...,k_n)\in \mathbb{Z}^n$, let
$v_k=v_{k_1}(1)\cdots v_{k_n}(n)$ where, for $1\leq i\leq n$ and
$m\geq 0$: $\,v_{m}(i)=x_i^m$, $\,v_{-m}(i)=y_i^m$, $v_0(i)=1$. It
follows from  the  definition of the  GWA that
$$A=\oplus_{k\in \mathbb{Z}^n} A_k$$
is a $\mathbb{Z}^n$-graded algebra ($A_kA_e\subseteq A_{k+e},$ for
all  $k,e \in \mathbb{Z}^n$),
 where $A_k=Dv_k=v_kD$.

 The tensor product (over the base field) $A\t A'$ of generalized Weyl algebras
 of degree $n$ and $n'$ respectively is a GWA of degree $n+n'$:
$$A\otimes A'=D\otimes D'((\sigma ,\sigma'), (a, a')).$$

Let $\CP_n$ be a polynomial algebra  $K[H_1, \ldots , H_n]$ in $n$
indeterminates  and let $\sigma=(\sigma_1,...,\sigma_n)$ be an
$n$-tuple  of  commuting automorphisms of $\CP_n$ such that
 $\s_i(H_i)=H_i-1$ and $\s_i(H_j)=H_j$, for $i\neq j$. Let $A_n=K\langle x_1, \ldots , x_n , \der_1, \ldots,
\der_n\rangle$ be the  Weyl algebra.  The algebra homomorphism
\begin{equation}\label{AnGWA}
A_n\ra \CP_n ((\sigma_1,...,\sigma_n), (H_1, \ldots , H_n)), \;\;
x_i\mapsto  x_i, \; \; \der_i\mapsto y_i, \; \; i=1, \ldots , n,
\end{equation}
  is an isomorphism. We identify the Weyl algebra $A_n$
with the GWA above via this isomorphism. Note that $H_i= \der_i
x_i = x_i\der_i+1$. Denote by $S_n$ the multiplicative submonoid
of $\CP_n$ generated by the elements $H_i+j$, $i=1, \ldots , n$,
and $j\in \mathbb{Z}$. It follows  from the above  presentation of
the Weyl algebra $A_n$ as a GWA that $S_n$ is an Ore set in $A_n$,
and, using the $\mathbb{Z}^n$-grading, that the (two-sided)
localization $\CA_n:=S_n^{-1}A_n$ of the Weyl algebra $A_n$ at
$S_n$ is the {\em skew Laurent polynomial ring}
\begin{equation}\label{Anskewlaurent}
\CA_n=S_n^{-1}\CP_n [x_1^{\pm 1}, \ldots ,x_n^{\pm 1};
\sigma_1,...,\sigma_n]
\end{equation}
with coefficients from
$$ S_n^{-1}\CP_n=K[H_1^{\pm 1}, (H_1 \pm 1)^{-1}, (H_1 \pm 2)^{-1}, \ldots ,
H_n^{\pm 1}, (H_n \pm 1)^{-1}, (H_n \pm 2)^{-1}, \ldots ],$$ the
localization of $\CP_n$ at $S_n$. We identify the Weyl algebra
$A_n$ with the subalgebra of $\CA_n$ via the monomorphism,
$$A_n\ra \CA_n, \;\; x_i\mapsto x_i,\;\; \der_i\mapsto  H_ix_i^{-1}, \;\;
i=1, \ldots , n.$$ Let $k_n$ be the $n$'th {\em Weyl skew field},
that is the full ring of quotients of the $n$'th Weyl algebra
$A_n$ (it exists by  Goldie's Theorem  since $A_n$ is a Noetherian
domain). Then the algebra $\CA_n$ is a $K$-subalgebra of $k_n$
generated by the elements $x_i$, $x_i^{-1}$, $H_i$ and $H_i^{-1}$,
$i=1, \ldots , n$ since, for all natural $j$,
$$ (H_i\mp j)^{-1}=x_i^{\pm j}H_i^{-1}x_i^{\mp j}, \;\; i=1, \ldots , n.$$
Clearly, $\CA_n \simeq \CA_1 \t \cdots \t \CA_1$ ($n$ times).

{\it Definition}. A $K$-algebra $R$ has the {\em endomorphism
property over} $K$ if, for each simple $R$-module $M$, $\End_R(M)$
is algebraic over $K$.

\begin{theorem}\label{bigAn}
\cite{Bav-Ann2001} Let $K$ be a field of characteristic zero.
\begin{enumerate}
\item The algebra $\CA_n$ is a simple, affine, Noetherian domain.
\item The Gelfand-Kirillov dimension $\GK (\CA_n)=3n$ $(\neq
2n=\GK (A_n))$. \item The (left and right) global dimension ${\rm
gl.dim} (\CA_n)=n$. \item The (left and right) Krull dimension
${\rm K.dim} (\CA_n)=n$. \item Let ${\rm d}={\rm gl.dim}$ or ${\rm
d}= {\rm K.dim}$. Let $R$  be a Noetherian $K$-algebra with ${\rm
d}(R)<\infty $ such that $R[t]$, the polynomial ring in a central
 indeterminate, has the endomorphism property over $K$. Then
${\rm d}(\CA_1\t R)= {\rm d}(R)+1$. If, in addition, the field $K$
is algebraically closed and uncountable, and the algebra $R$ is
affine, then
 ${\rm d}(\CA_n\t R)= {\rm d}(R)+n$.
\end{enumerate}
\end{theorem}

$\GK ({\cal A}_1)=3$ is due to A. Joseph \cite{Jos1}, p. 336;
 see also  \cite{KL}, Example 4.11, p. 45.

It is an experimental fact that many small quantum groups are
GWAs. More about GWAs and their generalizations the interested
reader can find in \cite{Alev-Far-La-Sol-00, Artam-Cohn-99,
Bav-inform, Bav-Jor-01, Bek-Ben-Fut-04,  Far-Sol-SuAl-03, Hart-06,
Kir-Mus-Pas-99, Kir-Kuz-05,  Maz-Tur-02, Prest-Pu-02, Rich-Sol-06,
Staf-sl2-82}.

{\bf Projections}. The polynomial algebra $P_n=\oplus_{\alpha \in
\N^n}Kx^\alpha$ is a left $\mA_n$-module and $\End_K(P_n)$-module.
 For each $i=1, \ldots , n$ and $\alpha \in \N^n$,
$H_i(x^\alpha )= (\alpha_i +1) x^\alpha$,  and so $H_i$ is an {\em
invertible} map with  $H_i^{-1} (x^\alpha )=(\alpha_i +1)^{-1}
x^\alpha$. Let $h_i:= x_i\der_i$. Then, in $P_n$, $h_i(x^\alpha )=
\alpha_i x^\alpha$, and so $\ker (h_i) = K[x_1, \ldots , \hx_i,
\ldots , x_n]$.

Note that (in $\End_K(P_n)$) 
\begin{equation}\label{ob1}
(H^{-1} \der )^i (xH^{-1})^i= \frac{1}{H(H+1)\cdots (H+i-1)}, \;\;
i\geq 1.
\end{equation}

For each $\alpha \in \N^n$, the following element of $\End_K(P_n
)$ is invertible,  
\begin{equation}\label{maa}
(-\alpha , \alpha ) := \der^\alpha x^\alpha = \prod_{i=1}^n H_i(
H_i+1) \cdots (H_i+\alpha_i-1).
\end{equation}

\begin{lemma}\label{p8Ap7}
Let $K$ be a commutative $\Q$-algebra and $\alpha \in \N_n$. Then
 $x^\alpha (-\alpha , \alpha )^{-1}
\der^\alpha$ is the projection onto the ideal $(x^\alpha )$ of
$P_n$ in the decomposition $P_n = (\oplus_{\beta : x^\beta \not
\in (x^\alpha )}Kx^\beta ) \oplus (x^\alpha )$.
\end{lemma}

{\it Proof}.  If $x^\beta \not \in (x^\alpha )$ then $\der^\alpha
( x^\beta ) =0$, and so $x^\alpha \frac{1}{\der^\alpha
x^\alpha}\der^\alpha (x^\beta ) =0$. If $x^\beta  \in (x^\alpha )$
then $x^\alpha \frac{1}{\der^\alpha x^\alpha}\der^\alpha (x^\beta
) = x^\alpha \frac{1}{\der^\alpha x^\alpha}\der^\alpha x^\alpha
(x^{\beta -\alpha })= x^\alpha (x^{\beta -\alpha })=x^\beta $.
$\Box $

Lemma \ref{p8Ap7} is useful in producing various projections onto
homogeneous $K$-submodules of $P_n$. Let $S$ be a subset of $\N^n$
and $S'$ be its complement. Then $P_n = P_{n,S}\oplus P_{n,S'}$
where $P_{n,S}:=\oplus_{\alpha \in S}Kx^\alpha$ and
$P_{n,S'}:=\oplus_{\alpha \in S'}Kx^\alpha$. Then $\pi_S:=
\sum_{\alpha \in S}\pi_\alpha$ is the projection onto $P_{n,S}$.

{\it Example}. Let $a, b\in \N^n$ with $a\leq b$,  i.e. $a_1\leq
b_1,\ldots , a_n\leq b_n$; and $C:= \{ \alpha \in \N^n\, | \,
a\leq \alpha \leq b\}$ be the discrete cube in $\N^n$ and $C'$ be
its complement. Then $\pi_C:= \sum_{\alpha \in C} \pi_\alpha$ is
the projection onto $P_{n, C}$ in the decomposition $P_n = P_{n,
C}\oplus P_{n,C'}$. Note that $\pi_C= \prod_{i=1}^n
(x_i^{a_i}\frac{1}{\der_i^{a_i}
x_i^{a_i}}\der_i^{a_i}-x_i^{b_i}\frac{1}{\der_i^{b_i}
x_i^{b_i}}\der_i^{b_i})$, by Lemma \ref{p8Ap7}. In more detail,
for each $i=1, \ldots , n$, $x_i^{a_i}\frac{1}{\der_i^{a_i}
x_i^{a_i}}\der_i^{a_i}$ is the projection onto the ideal
$x_i^{a_i}K[x_i]$ in the decomposition $K[x_i]=
(\oplus_{j=0}^{a_i-1} Kx_i^j)\oplus x_i^{a_i}K[x_i]$. Therefore,
$p_i:=x_i^{a_i}\frac{1}{\der_i^{a_i}
x_i^{a_i}}\der_i^{a_i}-x_i^{b_i}\frac{1}{\der_i^{b_i}
x_i^{b_i}}\der_i^{b_i}$ is the projection onto $Kx_i^{a_i}\oplus
Kx_i^{a_i+1} \oplus \cdots \oplus Kx_i^{b_i}$. Now, it is obvious
that the product $p_1\cdots p_n$ is equal to $\pi_C$.

{\bf The Jacobian algebra $\mA_n$}.  Let $K$ be a commutative
$\Q$-algebra.

{\it Definition}. The {\bf Jacobian  algebra} $\mA_n$ is the
subalgebra of $\End_K(P_n )$ generated by the Weyl algebra $A_n$
and the elements $H_1^{-1} , \ldots , H_n^{-1}$.

Surprisingly, the Weyl algebras $A_n$ and the Jacobian algebras
$\mA_n$ have little in common. For example,  the algebra $\mA_n$
contains the infinite direct sum $K^{(\N )}$ of rings $K$. In
particular, $\mA_n$ is not a domain, and we will see that $\mA_n$
is not left or right Noetherian algebra.

By the very definition, 
\begin{equation}\label{AAtn}
\mA_n = \mA_1(1) \t \mA_1(2) \t \cdots \t \mA_1(n) \simeq
\mA_1^{\t n},
\end{equation}
where $\mA_1(i) : = K\langle x_i, \der_i, H_i^{-1} \rangle$ and
$\t := \t_K$. The algebra $\mA_n$ contains all the integrations
$\int_i= x_iH_i^{-1}$, $1\leq i \leq n$. In the algebra $\mA_n$,
each element $\der_i$ has a right inverse,  $\int_i$:
$\der_i\int_i= {\rm id}_{P_n}$; and each element $x_i$ has a left
inverse, $H_i^{-1} \der_i$: $H_i^{-1}\der_i x_i={\rm id}_{P_n}$.
So, the algebra $\mA_n$ contains all necessary operations of
Analysis (like integrations and differentiations) to deal with
polynomials. The algebra $\mA_n$ contains all integro-differential
operators. By (\ref{AAtn}), properties of the algebra $\mA_n$ is
mainly determined by properties of the algebra $\mA_1$.

We pointed out already that the multiplicative submonoid $S_n$ of
$K[H_1, \ldots , H_n]$ generated by the elements $H_i+j$, $1\leq i
\leq n$, $j\in \Z$, is a (left and right) Ore set of the Weyl
algebra $A_n$ and $S_n^{-1} A_n= \CA_n$. Using the $\Z^n$-grading
of the Weyl algebra $A_n$ coming from its presentation as a
generalized Weyl algebra one can easily verify that the
multiplicative submonoid $S_{n,+}$ of $K[H_1, \ldots , H_n]$
generated by the elements $H_i+j$, $1\leq i \leq n$, $j\in \N$, is
not a (left and right) Ore set of the Weyl algebra $A_n$. This
also follows from the fact that the algebra $\CA_n$ is a domain
but $\mA_n$ is not (if $S_{n,+}$ were a left or right Ore set then
$\mA_n\subseteq \CA_n$, a contradiction).

Consider the case $n=1$. Let $2^\N$ be the Boolean algebra of all
subsets of $\N$ and $\CB_1$ be the Boolean subalgebra generated by
all the finite subsets of $\N$. So, a subset $S$ of $\N$ is an
element of $\CB_1$ iff either $S$ is finite or co-finite (that is,
 its complement is finite). Note that the Jacobian algebra $\mA_1$
 contains all the projections $\pi_S$, $S\in \CB_1$.

In order to make formulae more readable, we drop the subscript
$1$. So, let, for a moment, $x:=x_1$, $\der :=\der_1$, and $H:=
H_1$. Since $(H^{-1}\der )^i H^{-1} x^i= (H^{-1}\der )^i
x^i(H+i)^{-1}=(H+i)^{-1}$, $i\geq 1$, the algebra $\mA_1$ contains
the subalgebra $L:= K[ H, H^{-1}, (H+1)^{-1}, \ldots , (H+i)^{-1},
\ldots ]$. For each $i\geq 1$, $\frac{1}{\der^ix^i}=
\frac{1}{H(H+1)\cdots (H+i-1)}\in L$.  Let $\mD_1:=
L+\sum_{i,j\geq 1} Kx^iH^{-j}\der^i$ and $V:= \oplus_{j\geq 1}
KH^{-j}$. The next theorem gives a $K$-basis for the algebra
$\mA_1$.

\begin{theorem}\label{8Ap7}
Let $K$ be a commutative $\Q$-algebra. Then the Jacobian algebra
$\mA_1=\oplus_{i\in \Z} \mA_{1, i}$ is a $\Z$-graded algebra
($\mA_{1, i}\mA_{1, j}\subseteq \mA_{1, i+j}$ for all $i,j\in \Z$)
where $\mA_{1, 0}=\mD_1$,  $\mD_1=L\oplus (\oplus_{i,j\geq 1}
Kx^iH^{-j} \der^i)$; and, for each $i\geq 1$, $\mA_{1, i } = x^i
\mD_1$ and $\mA_{1,-i} = \mD_1\der^i$.
\end{theorem}

{\it Proof}. First, let us prove that the sum in the definition of
$\mD_1$ is the direct one. Suppose that $r:=
l+xv_1\der+x^2v_2\der^2+\cdots + x^sv_s\der^s=0$ is a nontrivial
relation for some elements $l\in L$ and $v_i\in V:= \oplus_{j\geq
1}KH^{-j}$. We seek a contradiction. Since $L$ is a subalgebra of
$\End_K(P_n )$, the
 relation $r$ is not of the type $r=l$. So, we can assume that
 $v_s\neq 0$ and the natural number $s\geq 1$ is called the degree
 of the relation $r$. Let $r$ be a nontrivial relation of the
 least degree. The rational function $l\in K(H)$ can be written as
 $\frac{p}{q}$ where $p$ and $q$ are co-prime polynomials and the
 polynomial $q$ is a finite product of the type $\prod_{i\geq
 0}(H+i)^{n_i}$. Evaluating the relation $r$ at $1$: $0= r(1) =
 l(1) = \frac{p(1)}{q(1)}$, we see that the polynomial $p$ is
 equal to $(H-1) p'$ for some polynomial $p'\in K[H]$. Suppose
 that $n=1$, then $0=H^{-1}\der r x= H\s^{-1} (\frac{p'}{q})+v_1H$
  where $\s : H\mapsto H-1$ is the
 $K$-automorphism of the polynomial algebra $K[H]$ (and of its field
 of fractions $K(H)$). It follows that $v_1=-\s^{-1}
 (\frac{p'}{q})\in V\cap \s^{-1} (L)=0$, a contradiction.
 Therefore, $s\geq 2$.

 The $K[H]$-module $(V+K[H])/K[H]\simeq V$ has the $K$-basis $\{
 H^{-i}, i\geq 1\}$. In this basis, the matrix of the $K$-linear
 map $v\mapsto (H+\l )v$ (where $\l \in K$) is an upper triangular
 infinite matrix with $\l $ on the diagonal. In particular, the
 matrices of the maps $H+1, H+2, \ldots , H+s-1$,  are invertible,
 upper triangular. It follows from this fact that the relation
 $$ H^{-1} \der r x = H\s^{-1} (\frac{p'}{q}) +v_1H +xv_2(H+1)
 \der +\cdots + x^{i-1} v_i (H+i-1) \der^{i-1} +\cdots + x^{s-1}
 v_s(H+s-1) \der^{s-1}$$
 has degree $s-1$ since $s\geq 2$. Then, by induction on $s$, and the
 fact that each matrix of the map $H+i-1$, $2\leq i \leq  s$, is
 an upper triangular, invertible matrix with $i-1\neq 0$ on the
 diagonal, we have $v_2=\cdots = v_s=0$ and $H\s^{-1} (\frac{p'}{q})
 +v_1H=0$. Then, $v_1= -\s^{-1} (\frac{p'}{q})\in V\cap \s^{-1}
 (L)=0$. This means that the relation $r$ is a trivial one, a
 contradiction. This finishes the proof of the claim.

Let $(G, +)$ be an additive group (not necessarily  commutative)
and $U= \oplus_{\alpha \in G}U_\alpha $ be  a $G$-graded
$K$-module, i.e. a direct sum of $K$-modules $U_\alpha$. A
$K$-linear map $f: U\ra U$ has degree $\beta \in G$ if $f(U_\alpha
) \subseteq U_{ \beta +\alpha} $ for all $\alpha \in G$. The set
$E_\beta$ of all $K$-linear maps of degree $\beta$ is a
$K$-submodule of $\End_K(U)$. Clearly, $E_\beta = \prod_{\beta \in
G} \Hom_K(U_\alpha , U_{\beta +\alpha} )$. In particular,
$E_0=\prod_{\alpha \in G}\End_K(U_\alpha )$.
 It follows at once that the sum $E:= \sum_{\beta \in
 G} E_\beta\subseteq \End_K(U)$ is a direct one, 
\begin{equation}\label{E=Eb}
E= \oplus_{\beta \in G} E_\beta\;\; {\rm and}\;\; E_\beta E_\g
\subseteq E_{\beta +\g }, \;\; \beta, \g \in G .
\end{equation}
So, $E:= E(U)$ is a $G$-graded ring.

 The $K$-module
$K[x]=\oplus_{i\geq 0}Kx^i$ is naturally $\Z$-graded (even
$\N$-graded). By (\ref{E=Eb}), the sum $S:= \sum_{i\geq 1}
\mD_1\der^i+\mD_1+\sum_{i\geq 1} x^i\mD_1$ is a direct sum since
the maps $\mD_1\der^i$, $\mD_1$, and $x^i\mD_i$ have degree $-i$,
$0$, and $i$ respectively. In order to prove that $\mA_1=
\oplus_{i\in \Z} \mA_{1,i}$ it suffices to show that
$\mA_1\subseteq S$. It follows directly from the inclusions:
\begin{eqnarray*}
 \mD_1x^i&\subseteq &x^i\mD_1, \;\;\;\;\;\;  \der^i \mD_1\subseteq \mD_1\der^i, \;\;\;\;\;\, i\geq 0,  \\
x^i\der^j&\subseteq &\mD_1\der^{j-i},\;\;\;\; x^j\der^i\subseteq
x^{j-i}\mD_1,\;\; j\geq i,
\end{eqnarray*}
that $\mA_1=\sum_{i,j\geq 0} x^i\mD_1\der^j=\sum_{s\in
\Z}(\sum_{i-j=s}x^i\mD_1\der^j)$. It remains to show that, for
$s\geq 0$, $\sum_{i-j=s}x^i\mD_1\der^j\subseteq x^s\mD_1$; and,
for $s<0$, $\sum_{i-j=s}x^i\mD_1\der^j\subseteq \mD_1\der^{-s}$.

Consider the case $s=0$. We have to show that $\sum_{i\geq
0}x^i\mD_1\der^i\subseteq \mD_1$. By the very definition of
$\mD_1$, this is equivalent to the inclusions $x^i
L\der^i\subseteq \mD_1$, $i\geq 0$; and, by the very definition of
$L$, this is equivalent to the inclusions $x^i(H+j)^{-k} \der^i\in
\mD_1$, $i,j, k\geq 0$.

If $i\leq j$ then $x^i(H+j)^{-k} \der^i= (H+j-i)^{-k} x^i\der^i\in
L\subseteq \mD_1$.

If $i> j$ then $x^i(H+j)^{-k} \der^i= x^{i-j}H^{-k}
x^j\der^i=x^{i-j} H^{-k} (H-1) (H-2) \cdots ( H-j) \der^{i-j} \in
x^{i-j}(V+K[H])\der^{i-j}\subseteq \mD_1$. This proves the case
$s=0$.

If $s\geq 1$ then $\sum_{i-j=s} x^i\mD_1\der^j = x^s(\sum_{k\geq
0} x^k\mD_1\der^k) \subseteq x^s\mD_1$, by the case $s=0$.

If $s\leq - 1$ then $\sum_{i-j=s} x^i\mD_1\der^j = (\sum_{k\geq 0}
x^k\mD_1\der^k) \der^{-s}\subseteq \mD_1\der^{-s}$, by the case
$s=0$.

Thus, the equality $\mA_1=\oplus_{i\in \Z} \mA_{1,i}$ is
established. By (\ref{E=Eb}), this is a $\Z$-graded algebra since
the maps $\mA_{1,i}$ have degree $i$.  $\Box $

Despite the fact that Theorem \ref{8Ap7} provides a cute $K$-basis
for the algebra $\mD_1$ it is unsuitable for computations: to
write down the product of the type $x^iH^{-j}\der^i\cdot
x^kH^{-l}\der^k$ literally takes half a page. Later, in Corollary
\ref{c8Ap7} a more conceptual $K$-basis is introduced, and which
is more important we interpret elements of $\mD_1$ as functions
from $\N$ to $K$. We will see that the ring $\mD_1$ is large and
has analytic flavour.

The polynomial algebra $K[x]= \oplus_{i\geq 0}Kx^i$ is naturally a
$\Z$-graded algebra. Let $E= \oplus_{i\in \Z}E_i$ be the algebra
from (\ref{E=Eb}) for $K[x]$. The $E$-module $K[x]$ is simple.
Note that the map
$$E_0= \{ f\in \End_K(K[x])\, | \, f(x^i)= f_ix^i, i\geq 0, f_i\in
K\}\ra  K^\N, \;\; f\mapsto (f_i),$$ is an  isomorphism of
$K$-algebras. In particular, $E_0$ is a commutative algebra, and
so $\mD_1$ is a {\em commutative} algebra since $\mD_1\subseteq
E_0$ ($K[H]$ is the faithful $\mA_1$-module). It is obvious that,
for $i\geq 0$, $E_i= x^iE_0$ and $E_{-i} = E_0\der^i$. The algebra
$K^\N$ is the algebra of all functions from $\N $ to $K$. When we
identify the set of monomials $\CM := \{ x^i\}_{i\in \N}$ and $\N$
via $x^i\mapsto i$ the algebras $E_0$ and $K^\N$ are identified.
So, each element of $E_0$ can be seen as a function. This is a
very nice observation indeed as we can use facts and terminology
of Analysis. For a function $\v:\N \ra K$, the set $\supp (\v
):=\{ i\in \N \, | \, \v (i)\neq 0\}$ is called the {\em support}
of $\v$. The set of functions $F_0$ with finite support is an
ideal of the algebra $E_0$. Clearly, $F_0= \oplus_{i\in \N}
K\pi_i$ where
$$\pi_i:=x^i\frac{1}{\der^ix^i}\der^i-x^{i+1}\frac{1}{\der^{i+1}x^{i+1}}\der^{i+1}:K[x]\ra K[x]$$
is the projection onto $Kx^i$ (Lemma \ref{p8Ap7}),  i.e.
$\pi_i(x^j) =\d_{ij} x^j$ where $\d_{ij}$ is the Kronecker delta.

Let $\pi_{-1}:=0$; then
$$ x\pi_i= \pi_{i+1}x, \;\; \der \pi_i= \pi_{i-1}\der, \;\; i\geq
0.$$

 The concept of support can be extended to an arbitrary element
 $f$ of the algebra $E$ as $\supp (f) := \{ i\in \N \, | \,
 f(x^i)\neq 0\}$. The set $F$ of all maps $f\in E$ with finite
 support is a $\Z$-graded algebra $F= \oplus_{i\in \Z}F_i$ without 1 where  $F_i=
 F\cap E_i$. For $i\geq 1$, $F_i= x^iF_0= \oplus_{j\in \N}Kx^i\pi_j=\oplus_{j\in \N}K\pi_{j+i}x^i=F_0x^i$ and
 $F_{-i}= F_0\der^i= \oplus_{j\in \N}K\pi_j\der^i=\oplus_{j\in
 \N}K\der^i\pi_{j+i}=\der^iF_0$. Note that $F= \{ f\in
 \End_K(K[x])\, | \, f(x^iK[x])=0$ for some $i\in \N\}$. It is
 obvious that $F$ is an {\em ideal} of the algebra $E$, and so $F$
 is also an ideal of $\mA_1$ since $F\subseteq \mA_1\subseteq E$.

 Note that $h:= x\der \in E_0$ and $h(x^i) = ix^i$, $i\geq 0$. So, $h$ can be
 identified with the function $\N \ra K$, $i\mapsto i$. Note that
 $H= h+1$. When $h$ runs through $0,1,2\ldots$, $H$ runs through
 $1,2, \ldots $. Under the identification $E_0= K^\N$, for
 $i,j\geq 1$,
\begin{equation}\label{xiHm1}
x^iH^{-j} \der^i = \begin{cases}
\frac{(H-1)(H-2)\cdots (H-i+1)}{(H-i)^{j-1}}& \text{if $H=i+1, i+2, \ldots $},\\
0& \text{if $H=1, 2, \ldots , i $}.\\
\end{cases}
\end{equation}
This means that the element $ x^iH^{-j} \der^i$ is a function of
the discrete argument $H=1,2, \ldots $ which takes zero value for
$H=1, 2, \ldots, i$; and $\frac{(H-1)(H-2)\cdots
(H-i+1)}{(H-i)^{j-1}}$ for  $H=i+1, i+2, \ldots $. Before the
identification this simply means that
$$x^iH^{-j} \der^i (x^k) = \begin{cases}
\frac{(k+1-1)(k+1-2)\cdots (k+1-i+1)}{(k+1-i)^{j-1}}x^k& \text{if $k\geq i$},\\
0& \text{if $k=0,1,  \ldots , i-1$}.\\
\end{cases}
$$
So, the function $x^iH^{-j} \der^i$ is almost a rational function.
 The case $i=j=1$ is rather special, it yields almost a constant
 function
$$xH^{-1} \der = \begin{cases}
1& \text{if $H=2, 3,  \ldots $},\\
0& \text{if $H= 1  $}.\\
\end{cases}
$$
Similarly, for each $i=1,2, \ldots$, the element $\rho_i:=
x^i\frac{1}{\der^ix^i}\der^i\in \mD_1$ is the function
\begin{equation}\label{xiHm2}
\rho_i= \begin{cases}
1& \text{if $h=i, i+1,  \ldots $},\\
0& \text{if $h=0, 1,  \ldots , i-1 $}.\\
\end{cases}
\end{equation}
For each $i\geq 0$,  $\pi_i= \rho_i-\rho_{i+1}\in\mD_1$ where
$\rho_0:=1$, and so $F_0= \oplus_{i\geq 0} K\pi_i\subseteq \mD_1$.
More generally, for each $i=1,2, \ldots $ and $j\in \N$, let
\begin{equation}\label{xiHm3}
\rho_{ji}:= x^i\frac{1}{H^j\der^ix^i}\der^i=
\frac{1}{(H-i)^j}\rho_i=\begin{cases}
\frac{1}{(H-i)^j}& \text{if $H=i+1, i+2, \ldots $},\\
0& \text{if $H=1, 2, \ldots , i $}.
\end{cases}
\end{equation}
Note that all $\rho_{ji}\in \mD_1$ and $\rho_{0i}=\rho_i$.  For
$\l \in \Z$, the element $H+\l $ is invertible in $\mD_1$ iff $\l
\neq -1, -2, \ldots $ iff $H+\l$ is invertible in $E_0$. For
$i=1,2,\dots $, $\ker_{\mD_1}(H-i)^j = \ker_{E_0}(H-i)^j = K\pi_i$
where $\ker_{\mD_1}(H-i)^j$ is the kernel of the map $\mD_1\ra
\mD_1$, $d\mapsto (H-i)^jd$. Similarly, $\ker_{E_0}(H-i)^j$ is
defined.

For each $i\geq 0$, let $\pi_i':= 1-\pi_i$. For natural numbers
$i,j\geq 1$, consider the element $\frac{1}{(H-i)^j}\pi_{i-1}'$ of
$E_0$ which is as a linear map defined by the rule
$$
\frac{1}{(H-i)^j}\pi_{i-1}'(x^k)= \begin{cases}
\frac{1}{(k+1-i)^j}x^k& \text{if $k\neq i-1$},\\
0& \text{if $k=i-1 $}.\\
\end{cases}
$$
As a function, it is almost the rational function
$\frac{1}{(H-i)^j}$ but at $H=i$ it takes value $0$ rather than
$\infty$ as the usual function $\frac{1}{(H-i)^j}$ does. All
$\frac{1}{(H-i)^j}\pi_{i-1}'\in \mD_1$ since
$$
\frac{1}{(H-i)^j}\pi_{i-1}'= \begin{cases}
\rho_{j1}& \text{if $i=1$},\\
\rho_{ji}+\sum_{k=0}^{i-2}\frac{1}{(k+1-i)^j}\pi_k& \text{if $i\geq 2$}.\\
\end{cases}
$$
For $i,j, n,m\geq 1$, $\frac{1}{(H-i)^n}\pi_{i-1}'\cdot
\frac{1}{(H-i)^m}\pi_{i-1}'= \frac{1}{(H-i)^{n+m}}\pi_{i-1}'$,
$(H-i)^m\cdot \frac{1}{(H-i)^n}\pi_{i-1}'=
\frac{1}{(H-i)^{n-m}}\pi_{i-1}'$; and for $j\geq 1$ such that
$j\neq i$
$$ \frac{1}{(H-i)^n}\pi_{i-1}'\cdot \frac{1}{(H-j)^m}\pi_{j-1}'
\subseteq \sum_{s=1}^nK\frac{1}{(H-i)^s}\pi_{i-1}'+
\sum_{t=1}^mK\frac{1}{(H-j)^t}\pi_{j-1}'+K\pi_{i-1}+K\pi_{j-1}.$$
Clearly, the set 
\begin{equation}\label{mL1def}
\mL_1:=L\oplus \oplus_{i,j\geq 1}K \frac{1}{(H-i)^j}\pi_{i-1}'
\end{equation}
is a $K$-submodule of $\mD_1$. $\mL_1$ is not an algebra, though
it is an algebra modulo $F_0$ which is isomorphic to the algebra $
K[H^{\pm 1}, (H\pm 1)^{-1}, (H\pm 2)^{-1}, \ldots , ]$ (see
Corollary \ref{c8Ap7} and (\ref{mDF})).

\begin{corollary}\label{c8Ap7}
Let $K$ be a commutative $\Q$-algebra,
$\rho_{ji}:=x^i\frac{1}{H^j\der^ix^i}\der^i$, $j\geq 0$, $i\geq
1$. Then $\mD_1 = \mL_1\oplus F_0= L\oplus (\oplus_{i\geq 1, j\geq
0}K\rho_{ji })$.
\end{corollary}

{\it Proof}. By Theorem \ref{8Ap7}, $\mD_1= L\oplus
(\oplus_{i,j\geq 1} Kx^iH^{-j}\der^i)$.  By (\ref{xiHm1}) and
$F_0\subseteq \mD_1$, we have $\mD_1\subseteq \mL_1+F_0$.
 By (\ref{xiHm3}) and
$F_0\subseteq \mD_1$, we have the opposite inclusion, and so
$\mD_1=\mL_1+F_0= \mL_1\oplus F_0$ since $\mL_1\cap F_0=0$.

The $K$-module $M:= L+\sum_{j\geq 0, i\geq 1} K\rho_{ji}$ contains
$F_0= \oplus_{i\geq 0}\pi_i$ since $\pi_i= \rho_i-\rho_{i+1}$ for
all $i\geq 0$ where $\rho_0:=1$; and $M= L'+F_0$ where $L':=
L+\sum_{j, i\geq 1} K\rho_{ji}$. Consider the factor module
$M/F_0$. By (\ref{xiHm3}), $\rho_{ji}\equiv
\frac{1}{(H-i)^j}\pi_{i-1}'\mod F_0$, hence $L'\equiv \mL_1\mod
F_0$. Since $\mD_1= \mL_1\oplus F_0$ and $M= L'+F_0$, we must have
the equality $\mD_1 = M$. To finish the proof of the second
equality of the corollary it suffices to show that $L'+F_0=
L\oplus (\oplus_{i,j\geq 1}K\rho_{ji})\oplus F_0$ since then the
quality $ M = L\oplus (\oplus_{i\geq 1,j\geq 0} K\rho_{ji})$
follows as
 $ F_0= \oplus_{i\geq 0 }K\pi_i$ and $\pi_i= \rho_i-\rho_{i+1}$.
 Let $l+\sum_{i, j\geq 1} \l_{ji} \rho_{ji}+f=0$ for some $l\in L$, $\l_{ji}\in
 K$, and $f\in F_0$. Taking this equality modulo $F_0$ yields
 $l=0$ and $\l_{ji}=0$ since $\rho_{ji} \equiv \frac{1}{(H-i)^j}\pi_{i-1}'\mod
 F_0$. This implies $f=0$, and we are done.  $\Box $

Note that $F_0$ is an ideal of $\mD_1$ such that $F_0^2=F_0$ and
\begin{equation}\label{mDF}
\mD_1/F_0\simeq K[H^{\pm 1}, (H\pm 1)^{-1}, (H\pm 2)^{-1}, \ldots
].
\end{equation}
The equality $\mD_1= \mL_1\oplus F_0$ (Corollary \ref{c8Ap7})
means that the set
$$ \{ H^i, \frac{1}{(H+j)^k}, \frac{1}{(H-j)^k}\pi_{j-1}', \pi_l\, | \,
i\in \Z, l\in \N, j,k\geq 1\}$$ is a $K$-basis for $\mD_1$.
Clearly, the set 
\begin{equation}\label{mDF1}
 \{ H^i, \frac{1}{(H+j)^k}, \frac{1}{(H-j)^k}\pi_{j-1}'\, | \,
i\in \Z,  j,k\geq 1\}
\end{equation}
 is a $K$-basis for $\mL_1$. Similarly, the
equality $\mD_1= L\oplus (\oplus_{i\geq 1, j\geq 0}K\rho_{ji})$
means that the set
$$ \{ H^j, \frac{1}{(H+j)^k}, \rho_{ji}\, | \,
j\in \N,  i,k\geq 1\}$$ is a $K$-basis for $\mD_1$.

Clearly, $F=\oplus_{i,j\geq 0}KE_{ij}$ where $E_{ij}(x^k):=
\d_{jk} x^i$, i.e. $\{ E_{ij}\}$ are the `elementary matrices'
($E_{ij}E_{kl}=\d_{jk} E_{il}$), and 
\begin{equation}\label{Eijxd}
E_{ij} = \begin{cases}
x^{i-j}\pi_j& \text{if $i\geq j$},\\
(\frac{1}{H}\der )^{j-i} \pi_j& \text{if $i<j$}.\\
\end{cases}
\end{equation}
For $k\in \N$, $F_{\pm k}= \oplus_{i-j=\pm k} K E_{ij}$.
 Note that $E_{ij} = E_{ik}\pi_kE_{kj}$ for all $i,j,k\geq 0$.
Therefore, $F$ is a {\em simple} ring such that $F^2=F$, and
$K[x]$ is a simple faithful $F$-module. The ring $F$ is neither
left nor right Noetherian  as the next arguments show: for each
natural $k\geq 0$, let $L_k:=\oplus_{i\in \N, 0\leq j \leq
k}KE_{ij}$ and $R_k:= \oplus_{0\leq i \leq k, j\in \N} KE_{ij}$
then $L_0\subset L_1\subset \cdots$ and $R_0\subset R_1\subset
\cdots$ are strictly ascending sequences of left and right
$F$-modules respectively. This is the main reason why the Jacobian
algebra $\mA_1$ is also neither left nor right Noetherian (Theorem
\ref{12Ma7}.(3)).  The ring $F$ is neither left nor right
Artinian: for each natural $k\geq 1$, let $L_k':=\oplus_{i, j\in
\N}KE_{i,j2^k}$ and $R_k:= \oplus_{i, j\in \N} KE_{i2^k, j}$ then
$L_0'\supset L_1'\supset \cdots$ and $R_0'\supset R_1'\supset
\cdots$ are strictly descending sequences of left and right
$F$-modules respectively.

\begin{theorem}\label{12Ma7}
Let $K$ be a commutative $\Q$-algebra. Then
\begin{enumerate}
\item $\mA_1= \oplus_{i\geq 1} \mL_1\der^i\oplus \mL_1\oplus
(\oplus_{i\geq 1} x^i\mL_1)\oplus F$. \item The set $\{ H^i\der^l,
\frac{1}{(H+j)^k}\der^l, \frac{1}{(H-j)^k}\pi_{j-1}'\der^l,
x^mH^i, x^m\frac{1}{(H+j)^k}, x^m\frac{1}{(H-j)^k}\pi_{j-1}',
E_{st}\, | \, i\in \Z ; j,k,l\geq 1; m, s,t\in \N\}$ is a
$K$-basis for $\mA_1$. \item The algebra $\mA_1$ is neither a left
nor right Noetherian algebra.\item $F$ is an ideal of $\mA_1$,
$F^2=F$, and the factor algebra $\mA_1/F$ is canonically
isomorphic to the algebra $\CA_1$ (the localization of the Weyl
algebra $A_1$ at $S_1$, the multiplicative monoid generated by
$H+i$, $i\in \Z$).
\end{enumerate}
\end{theorem}

{\it Proof}. 1. By Corollary \ref{c8Ap7}, $\mD_1 = \mL_1\oplus
F_0$. For each natural number $i\geq 1$, $\mD_1\der^i =
\mL_1\der^i \oplus F_{-i}$ and $x^i\mD_1= x^i\mL_1\oplus F_i$.
Using these equalities together with the equalities $\mA_1=
\oplus_{i\geq 1}\mD_1\der^i\oplus \mD_1\oplus (\oplus_{i\geq 1}
x^i\mD_1)$ (Theorem \ref{8Ap7}) and $F= \oplus_{i\geq
1}F_0\der^i\oplus F_0\oplus (\oplus_{i\geq 1} x^iF_0)$, one
obtains the equality of statement 1 and (\ref{mDF1}).

2. Since, for all $i\geq 1$, the maps $\mL_1\ra \mL_1\der^i$,
$u\mapsto u \der^i$, and $\mL_1\ra x^i\mL_1$, $u\mapsto x^i u $,
are isomorphisms of $K$-modules, statement 2 follows from
statement 1.

3. For each $i\in \N$, the sum $I_i:= \oplus_{j=0}^iK\pi_j$ is an
ideal of $E_0$. Since $I_i\subseteq \mD_1\subseteq E_0$, one has
the strictly ascending chain of ideals of $\mD_1$: $I_0\subset
I_1\subset \cdots$. The ascending chain $\mA_1I_0\subset
\mA_1I_1\subset \cdots$ of left homogeneous ideals of the algebra
$\mA_1$ is strictly ascending since the zero component of the left
ideal $\mA_1 I_j= \oplus_{i\geq 1} \mD_1 \der^iI_j\oplus I_j\oplus
(\oplus_{i\geq 1} x^i I_j)$ is $I_j$ (note that
$\mD_1\der^iI_j\subseteq F_{-i}$ and $x^iI_j\subseteq F_i$).
Therefore, $\mA_1$ is not a left Noetherian algebra.

Similarly, the  ascending chain $I_0\mA_1\subset I_1\mA_1\subset
\cdots$ of right homogeneous ideals of the algebra $\mA_1$ is
strictly ascending since the zero component of the right ideal
$I_j\mA_1 = \oplus_{i\geq 1}  I_j\der^i\oplus I_j\oplus
(\oplus_{i\geq 1} I_jx^i \mD_1)$ is $I_j$, and so  $\mA_1$ is not
a right Noetherian algebra.

4. We proved already that $F$ is an ideal of $\mA_1$ such that
$F^2=F$. The $K$-module $K[x]$ is a topological $K$-module (even a
topological  $K$-algebra) with respect to the $\gm$-adic topology
determined by the $\gm$-adic filtration $\{ \gm^i\}_{i\geq 0}$ on
$K[x]$ where $\gm := (x)$. Let $\End_{K,c}(K[x])$ be the algebra
of all {\em continuous} $K$-endomorphisms of $K[x]$. Then
$E\subseteq\End_{K,c}(K[x])$. Let $\CG$ be the algebra of {\em
germs} of continuous $K$-endomorphisms of $K[x]$ at $0$. An
element of $\CG$ is an equivalence class $[f]$ of a continuous
$K$-linear map of the type $f: \gm^i\ra K[x]$, and two such maps
are equivalent, $f\sim f'$,  if they have the same restriction
$f|_{\gm^j} = f'|_{\gm^j}$ for a sufficiently large $j$ ($\gm^i$
is a topological $K$-module with respect to the induced topology
coming from the inclusion $\gm^i\subseteq K[x]$).

The kernel of the $K$-algebra homomorphism $E\ra \CG$, $ f\mapsto
[f]$, is $F$. Thus, the kernel of the $K$-algebra homomorphism $
g: \mA_1\ra \CG$, $f\mapsto [f]$, is also $F$ since $F\subseteq
\mA_1$. The image $g(\mD_1)$ is naturally isomorphic to the
algebra $S_1^{-1} K[H]$ since $g(H) = [H]$ and $g(
\frac{1}{(H-j)^k}\pi_{j-1}') = [\frac{1}{(H-j)^k}]$ for all
$j,k\geq 1$. Now, it follows from statement 1 and the
decomposition $\CA_1= \oplus_{i\geq 1}S_1^{-1} K[H]\der^i\oplus
S_1^{-1} K[H]\oplus (\oplus_{i\geq 1} x^iS_1^{-1} K[H])$ that the
image $g(\mA_1)$ is naturally isomorphic to the algebra $\CA_1$.
$\Box $

An ideal $I$ of a ring $R$ such that $0\neq I\neq R$ is called a
{\em proper} ideal of $R$.
\begin{corollary}\label{a12Ma7}
 Let $K$ be a
 field of characteristic zero. Then
\begin{enumerate}
\item $F$ is the only proper ideal of the algebra $\mA_1$, hence
$F$ is a maximal ideal. \item $\mA_1/F\simeq \CA_1$ is a simple
Noetherian domain. \item $\GK ( \mA_1) = \GK (\CA_1) = 3$. \item
$\mA_1$ is a prime ring. \item The algebra $\mA_1$ is central,
i.e. the centre of $\mA_1$ is $K$.
\end{enumerate}
\end{corollary}

{\it Proof}. 2. Statement 2 follows from Theorem \ref{12Ma7}.(4)
and Theorem \ref{bigAn}.(1).

1. By statement 2, $F$ is a maximal ideal of $\mA_1$. Let $I$ be a
proper ideal of $\mA_1$. We have to show that $I=F$. Let $a$ be a
nonzero element of $I$. Then $0\neq FaF\subseteq F\cap I$, and so
$FaF = F$ since $F$ is a simple algebra (i.e. a simple
$F$-bimodule). Now, $F\subseteq I$ implies $F=I$ by the maximality
of $F$. So, $F$ is the only proper ideal of the algebra $\mA_1$.

3. Since $\mA_1/F\simeq \CA_1$, we have $\GK (\mA_1)\geq \GK
(\CA_1) = 3$ (Theorem \ref{bigAn}.(2)). Since $\mA_1=
\oplus_{i\geq 1} \mL_1\der^i\oplus \mL_1\oplus (\oplus_{i\geq 1}
x^i\mL_1) \oplus F$ and $\CA_1= \oplus_{i\geq 1}
S_1^{-1}K[H]\der^i\oplus S_1^{-1}K[H] \oplus (\oplus_{i\geq 1}
x^iS_1^{-1}K[H])$, the reverse inequality $\GK (\mA_1)\leq 3$
follows from Theorem \ref{12Ma7}.(1,2) using the same sort of
estimates as in the proof of the inequality $\GK(\CA_1)\leq 3$
(see \cite{Bav-Ann2001} for details and the fact that the elements
of $F$ do not contribute to the growth of degree 3).

4. $F$ is the only proper ideal of $\mA_1$, and so $F^2=F$, hence
$\mA_1$ is a prime ring.

5. The field $K$ belongs to the centre of $\mA_1$. Let $z$ be a
central element of $\mA_1$. We have to show that $z\in K$. The
algebra $\mA_1/F\simeq \CA_1$ is central, hence $z=\l +f$ for some
$\l \in K$ and $f\in F$. Then $ z-\l = f$ belongs to the centre of
$F$ which is obviously  equal to zero. Hence $z= \l \in K$, as
required.  $\Box$

For $k\geq 1$,
\begin{eqnarray*}
 x\frac{1}{(H-j)^k}\pi_{j-1}'&=& \frac{1}{(H-1-j)^k}\pi_{j}'x, \;\; j\geq 1,  \\
 \der \frac{1}{(H-j)^k}\pi_{j-1}'&=& \frac{1}{(H+1-j)^k}\pi_{j-2}'\der, \;\; j\geq 2,  \\
 \der \frac{1}{(H-1)^k}\pi_0'&=& \frac{1}{H^k}\der .
\end{eqnarray*}

By (\ref{AAtn}), the algebra $\mA_n$ is the tensor product of the
$\Z$-graded algebras $\mA_1(i)=\oplus_{j\in \Z}\mA_{1,j}(i)$.
Therefore, the algebra $\mA_n$ is a $\Z^n$-graded algebra,
$$ \mA_n= \oplus_{\alpha \in \Z^n}\mA_{n, \alpha}, \;\;
\mA_{n,\alpha } := \otimes_{i=1}^n \mA_{1,\alpha_i}(i).$$ The
$\Z^n$-grading on $\mA_n$ is the tensor product of the
$\Z$-gradings of the tensor multiples, and an element $a$ of
$\mA_n$ belongs to $\mA_{n,\alpha }$ iff $a(x^\beta ) \in
Kx^{\alpha +\beta}$ for all $\beta \in \N^n$. Let
$$ \mD_n:= \mA_{n,0}= \mD_1(1)\otimes \mD_1(2)\otimes \cdots
\otimes \mD_1(n).$$ The polynomial algebra $P_n= \oplus_{\alpha
\in \N^n}P_{n,\alpha}$ is an $\N^n$-graded, hence a $\Z^n$-graded
algebra. Let $E= E(P_n)=\oplus_{\alpha \in \Z^n}E_\alpha$ be the
$\Z^n$-graded algebra as in (\ref{E=Eb}). The map
$$E_0=\{ f\in \End_K(P_n)\, |\, f(x^\alpha )= f_\alpha x^\alpha,  f_\alpha \in K, \alpha \in \N^n\}
 \ra K^{\N^n}, \;\; f\mapsto (f_\alpha ),$$ is
a $K$-algebra isomorphism. Each element $\alpha =
\sum_{i=1}^n\alpha_ie_i\in \Z^n= \oplus_{i=1}^n\Z e_i$ is a unique
difference  $\alpha = \alpha_+-\alpha_-$ where
$\alpha_+=\sum_{\alpha_i\geq 0}\alpha_ie_i$ and
$\alpha_-=-\sum_{\alpha_i\leq 0}\alpha_ie_i$. For each $\alpha \in
\Z^n$, $E_\alpha = x^{\alpha_+}E_0\der^{\alpha_-}$, and so
\begin{equation}\label{E=nEa}
E=\oplus_{\alpha \in \Z_n}x^{\alpha_+}E_0\der^{\alpha_-}.
\end{equation}
For each $\alpha \in \Z^n$, $\mA_{n,\alpha}= \mA_n\cap E_\alpha =
x^{\alpha_+}\mD_n\der^{\alpha_-}$ where $\mD_n= \mA_{n,0}=
\mA_n\cap E_0$.

{\bf The involution $\th$ on $\mA_n$}. Let $K$ be a commutative
$\Q$-algebra. The Weyl algebra $A_n$ admits the {\em involution}
$$\th : A_n\ra A_n, \;\; x_i\mapsto \der_i, \;\; \der_i\mapsto
x_i, \;\; i=1, \ldots , n,$$ i.e. it is a $K$-algebra
anti-isomorphism ($\th (ab) = \th (b) \th (a)$) such that  $\th^2=
{\rm id}_{A_n}$. The involution $\th$ can be uniquely extended to
the involution of $\mA_n$ by the rule 
\begin{equation}\label{thinv}
\th : \mA_n\ra \mA_n, \;\; x_i\mapsto \der_i, \;\; \der_i\mapsto
x_i,\;\;  \th (H_i^{-1})= H_i^{-1} \;\; i=1, \ldots , n.
\end{equation}
Uniqueness is obvious: $\th (H_i) = \th (\der_ix_i)= \th (x_i )
\th (\der_i) = \der_ix_i= H_i$ and so $\th (H_i^{-1}) = H_i^{-1}$.
To prove existence recall that each right module over a ring $R$
is a left module over the opposite ring $R^{op}$. The involution
$\th$ on $A_n$ comes from considering the polynomial algebra $P_n$
as the {\em right} $A_n$-module by the rule $pa:= \th (a) p$ for
all $p\in P_n$ and $a\in A_n$. Since $\th (H_i) = H_i$, $i=1,
\ldots , n$, $P_n$ is the {\em faithful} right $\mA_n$-module, and
this proves the existence of the involution $\th :\mA_n\ra \mA_n$
($\th$ is injective since $P_n$ is a faithful right
$\mA_n$-module, $\th$ is obviously surjective). So, the algebra
$\mA_n$ is {\em self-dual} (i.e. it is isomorphic to its opposite
algebra, $\th : \mA_n\simeq \mA_n^{op}$). This means that  left
and right algebraic  properties of the algebra $\mA_n$ are the
same.

For $n=1$, $F$ is the only proper ideal of $\mA_1$, hence $\th (F)
= F$. Moreover, 
\begin{equation}\label{thEij}
\th (E_{ij})=\frac{i!}{j!}\, E_{ji}
\end{equation}
where $0!:=1$.  In more detail, since $ \th (H)=H$ and  $E_{ii} =
\pi_i = x^i\frac{1}{(-i,i)} \der^i-x^{i+1}\frac{1}{(-i-1,i+1)}
\der^{i+1}$, we have $\th (E_{ii}) = E_{ii}$. For $ i>j$, $E_{ij}=
x^{i-j} \pi_j$, and so $$\th ( E_{ij} )= \pi_j\der^{i-j}=i(i-1)
\cdots (j+1)E_{ji}=\frac{i!}{j!}E_{ji}.$$ For $i<j$,
$E_{ij}=(\frac{1}{H}\der )^{j-i} \pi_j$, and so $\th (E_{ij})=
\pi_j(x\frac{1}{H})^{j-i}=\frac{1}{(i+1)(i+2)\cdots j}
E_{ji}=\frac{i!}{j!}E_{ji}$.

 For $n=1$, the ring $F=\oplus_{i,j\in \N} KE_{ij}$
is equal to the matrix ring $M_\infty (K):= \cup_{d\geq 1} M_d(K)$
where $M_d(K):=\oplus_{0\leq i,j\leq d-1}KE_{ij}$. The ring
$F=M_\infty (K)$ admits the canonical involution which is the
 transposition $(\cdot )^t: E_{ij}\mapsto E_{ji}$. Let $D_!$ be the infinite
 diagonal matrix ${\rm diag} (0!, 1!, 2!, \ldots )$. Then, for $u\in
 F=M_\infty (K)$,
\begin{equation}\label{1thEij}
\th (u)=D_!^{-1} u^tD_!.
\end{equation}
Note that $D_!\not\in M_\infty (K)$.

For an arbitrary $n$, $F^{\t n } = \oplus_{\alpha , \beta \in
\N^n} K E_{\alpha \beta}=M_\infty (K)^{\t n }$ where $E_{\alpha
\beta}:= \t_{i=1}^n E_{\alpha_i\beta_i}$. By (\ref{thEij}),
\begin{equation}\label{nthEij}
\th (E_{\alpha \beta})= \frac{\alpha !}{\beta !}\, E_{\beta
\alpha},
\end{equation}
\begin{equation}\label{thFn}
\th (F^{\t n } ) =  F^{\t n }.
\end{equation}
Let $D_{n,!}:= D_!^{\t n }$. Then, for $u\in F^{\t n }$,
\begin{equation}\label{2thFn}
\th (u) =D_{n,!}^{-1}u^tD_{n,!}
\end{equation}
where $(\cdot )^t:M_\infty (K)^{\t n } \ra M_\infty (K)^{\t n }$,
$E_{\alpha \beta}\mapsto E_{\beta \alpha}$, is the transposition
map.

Consider the bilinear, symmetric, non-degenerate form $(\cdot,
\cdot ): P_n\times P_n\ra K$ given by the rule $(x^\alpha ,
x^\beta ) := \alpha ! \d_{\alpha, \beta }$ for all $\alpha , \beta
\in \N^n$. Then, for all $p,q\in P_n$ and $a\in \mA_n$,
\begin{equation}\label{paq=t}
(p,aq)= (\th (a) p , q).
\end{equation}

The Weyl algebra $A_n$ admits, so-called, the {\em Fourier
transform}, it is the $K$-algebra automorphism $\CF : A_n\ra A_n$,
$x_i\mapsto \der_i$, $\der_i\mapsto - x_i$, $i=1, \ldots , n$.
Since $\CF (H_i) = -(H_i-1)$, $H_i$ is a unit of $\mA_n$ and
$H_i-1$ is not,  one {\em cannot} extend the Fourier transform to
$\mA_n$.

{\bf The algebra $\mA_n$ is a prime algebra}.  Consider the ideals
of the algebra $\mA_n$:
$$\gp_1:=F\t \mA_{n-1}, \gp_2:= \mA_1\t F\t \mA_{n-2}, \ldots ,
\gp_n:= \mA_{n-1} \t F.$$ Then $\mA_n/\gp_i\simeq (\mA_1/F)\t
\mA_{n-1} \simeq \CA_1\t \mA_{n-1}$ and $\cap_{i=1}^n \gp_i =
F^{\t n }$. Let $\ga_n:= \gp_1+\cdots +\gp_n$. Then
\begin{equation}\label{AnAn}
\mA_n/ \ga_n\simeq (\mA_1/F)^{\t n } \simeq \CA_1^{\t n } = \CA_n.
\end{equation}

\begin{corollary}\label{15Ma7}
Let $K$ be a field of characteristic zero. Then
\begin{enumerate}
\item $\GK (\mA_n) = 3n$. \item $\GK (M)\geq n $ for all nonzero
finitely generated (left or right) $\mA_n$-modules $M$.\item The
centre of $\mA_n$ is $K$. \item $F^{\t n } $ is the smallest
nonzero ideal of the algebra $\mA_n$, $(F^{\t n })^2= F^{\t n
}$.\item The algebra $\mA_n$ is prime. \item If $0\neq a\in \mA_n$
and $I$ is a nonzero ideal of $\mA_n$ then $aI\neq 0$, $Ia\neq 0$,
and $IaI\neq0$. \item The ideal $\ga_n$ is the largest  ideal of
$\mA_n$ distinct from $\mA_n$; $\ga_n^2=\ga_n$;  hence $\ga_n$ is
the only maximal ideal of $\mA_n$. \item ${}_{\mA_n}F^{\t n }
\simeq P_n^{(\N^n)}$ is a faithful, semi-simple, left
$\mA_n$-module; $F^{\t n }_{\mA_n} \simeq {P_n^{(\N^n)}}_{\mA_n}$
is a faithful, semi-simple, right $\mA_n$-module; ${}_{\mA_n}F^{\t
n }_{\mA_n} \simeq P_n^{(\N^n)}$ is a faithful, simple
$\mA_n$-bimodule. \item $F^{\t n}$ is the socle of $\mA_n$
considered as a left $\mA_n$-module, or a right $\mA_n$-module, or
an $\mA_n$-bimodule. \item ${}_{\mA_n}P_n$ (resp. $(P_n)_{\mA_n}$)
is the only faithful, simple, left (resp. right) $\mA_n$-module.
\item $\GK (\mA_n/ \ga ) = 3n$ for all ideals $\ga$ of $\mA_n$
such that $\ga \neq \mA_n$.
\end{enumerate}
\end{corollary}

{\it Proof}. 1. On the one hand, $\GK (\mA_n) \geq \GK (\mA_n/
\ga_n) = \GK (\CA_n ) = 3n$ (Theorem \ref{bigAn}.(2)); on the
other, $\GK (\mA_n) = \GK (\mA_1^{\t n }) \leq n \, \GK (\mA_1)=
3n$. Therefore, $\GK (\mA_n) = 3n$.

2. Statement 2  is an easy corollary of the {\em inequality of
Bernstein}: $\GK (N)\geq n $ for all nonzero finitely generated
(left or right) $A_n$-modules $N$. Let $M$ be a nonzero finitely
generated $\mA_n$-module and $0\neq u\in M$. Then $ \GK_{\mA_n} (
M) \geq  \GK_{\mA_n} (\mA_n u ) \geq  \GK_{A_n} (A_nu) \geq n$.

3. To prove that the centre $Z(\mA_n) $ of $\mA_n$ is $K$ we use
induction on $n$. The case $n=1$ is Corollary \ref{a12Ma7}.(5).
Suppose that $n>1$ and the algebra $\mA_m$ is central for all
$m<n$. The kernel of the algebra homomorphism $ \mA_n\ra
\prod_{i=1}^n \mA_n/ \gp_i$ is $\cap_{i=1}^n \gp_i = F^{\t n }$.
Since all the algebras $ \mA_n/\gp_i\simeq \CA_1\t \mA_{n-1}$ are
central, we have $Z(\prod_{i=1}^n \mA_n/ \gp_i) = \prod_{i=1}^n
Z(\mA_n/ \gp_i) = \prod_{i=1}^n K$ where $Z(R)$ is the centre of
$R$. If $z\in Z(\mA_n)$ then $z+F^{\t n } \in Z(\prod_{i=1}^n
\mA_n/ \gp_i) =\prod_{i=1}^n K$, and so $z= \l +f$ for some $\l
\in K$ and $f\in F^{\t n }$. Now, $f= z-\l \in Z(F^{\t n })=0$,
i.e. $z=\l \in K$, and so the algebra $\mA_n$ is central.

4. Clearly, $(F^{\t n })^2= (F^2)^{\t n }= F^{\t n }$. It remains
to prove  minimality of $F^{\t n }$. This is obvious for $n=1$
(Corollary \ref{a12Ma7}.(1)). To prove the general case we use
induction on $n$. Suppose that $n>1$ and the result is true for
all $n'<n$. Let $I$ be a nonzero ideal of $\mA_n$. We have to show
that $F^{\t n } \subseteq I$. Choose a nonzero element, say $a$,
of $I$. Since $ a\in \mA_n = \mA_1\t \mA_{n-1}$, the element $a$
can be written as a sum $a= \sum_{i=1}^s a_i\t b_i$ for some
elements $a_i\in \mA_1$ and $b_i\in \mA_{n-1}$ such that the
elements $a_1, \ldots , a_s$ and $b_1, \ldots  , b_s$ are
$K$-linearly independent elements of the algebras $\mA_1$ and
$\mA_{n-1}$ respectively. Choose an element, say $f$ of $F$ such
that $fa_1\neq 0$. Changing $a$ for $fa\neq 0$ (and deleting zero
terms of the type $fa_i\t b_i$) one may assume that all $a_i\in
F$. Note that $\{ E_{kl}\}$ is the $K$-basis of $F$. So, the
element $a$ can be written as  a finite sum $a= E_{kl}\t
\alpha+E_{st}\t \beta +\cdots +E_{pq}\t \g$ for some $K$-linearly
independent elements $\alpha , \beta, \ldots , \g$ of $\mA_{n-1}$
and distinct elements $E_{kl}, E_{st}, \ldots , E_{pq}$. Then $b:=
E_{kk}aE_{ll}= E_{kl} \t \alpha \in I$. By induction, $F^{\t
(n-1)}\subseteq \mA_{n-1} \alpha \mA_{n-1}$, and so
$$ I\supseteq \mA_n b \mA_n= \mA_1 E_{kl}\mA_1\t \mA_{n-1} \alpha
\mA_{n-1}\supseteq F\t F^{\t (n-1)}= F^{\t n }.$$
5. Let $I$ and
$J$ be nonzero ideals of $\mA_n$. By statement 4, they contain the
ideal $F^{\t n }$. Now, $IJ\supseteq (F^{\t n })^2= F^{\t n } \neq
0$. This means that $\mA_n$ is a prime ring.

6. Statement 6 follows directly from statement 5. Suppose that
$aI=0$ for some nonzero element $a $ of $\mA_n$ and a nonzero
ideal $I$. We seek a contradiction. Then $ 0 =\mA_naI= \mA_n a
\mA_nI\neq 0$ since $\mA_n$ is a  prime algebra, a contradiction.
 Similarly, $Ia=0$ (resp. $IaI=0$) implies $0= I\mA_n a\mA_n \neq 0$
 (resp. $0= I\mA_n aI\neq 0$), a contradiction.

7. $\ga_n\supseteq \ga_n^2= (\gp_1+\cdots +\gp_n)^2\supseteq
\gp_1^2+\cdots +\gp_n^2= \gp_1+\cdots +\gp_n = \ga_n$ since
$\gp_1^2=\gp_1, \ldots , \gp_n^2=\gp_n$, and so $\ga_n = \ga_n^2$.

It remains to show that $\ga_n$ is the largest ideal distinct from
$\mA_n$, that is $a\not\in \ga_n$ implies $\mA_n a\mA_n = \mA_n$
where $a\in \mA_n$. Let $B_n$ be the $K$-basis for $\mA_n$ that is
the tensor product of the $K$-bases from Theorem \ref{12Ma7}.(2).
 $B_n$ is the disjoint union of its two subsets $\CM_n:= \{ b\in
 B_n \, | \, b\in \ga_n\} $ and $\CN_n :=\{ b\in B_n \, | \,
 b\not\in \ga_n\}$.  So, $b\in \CM_n$ iff the product $b$ contains
  a matrix unit $E_{st}(i)\in F(i)$.  Clearly, $\ga_n =
  \oplus_{\mu \in \CN_n} K\nu$ and $\mA_n = \ga_n'\oplus \ga_n$
  where $\ga_n' := \oplus_{\mu \in \CM_n} K\mu$. Elements of the
  set $\mA_n \backslash \ga_n$ are called {\em generic}. So, an
  element of $\mA_n$ is generic iff it has at least one nonzero
  $\mu$-coordinate for some $\mu \in \CM_n$. We have to show that
  $\mA_n a\mA_n= \mA_n$ for all generic elements $a\in \mA_n$.
This is true when $n=1$ (Corollary \ref{a12Ma7}.(1)). To prove
general case we use induction on $n$. So, let $n\geq 2$ and we
assume that the claim is true for all $n'<n$. Let $a$ be a generic
element of $\mA_n$. Then $a= a_1\t b_1+\cdots +a_s\t b_s\in \mA_1
\t \mA_{n-1}$ where $a_i$ are nonzero elements of $\mA_1$ such
that $a_1$ is generic; $b_i$ are distinct elements of the basis
$B_{n-1}$ such that $b_1\in \CN_{n-1}$. By statement 6, $Fa_1F\neq
0$, and so $E_{ij} a_1E_{kl}\neq 0$ for some $i,j,k,l\in \N$.
Then, for each $t=1, \ldots , s$, $E_{ij} a_tE_{kl}= \l_t E_{il}$
for some $\l_t\in K$, necessarily $\l_1\neq 0$. Now, $E_{ij}
aE_{kl} = E_{il} \t u$ where $u:= \l_1b_1+\cdots +\l_sb_s$  is a
generic element of $\mA_{n-1}$ since $\l_1\neq 0$ and $b_1\in
\CN_{n-1}$. By induction, $\mA_{n-1} u \mA_{n-1} = \mA_{n-1}$, and
so
$$ \mA_na\mA_n \supseteq \mA_n ( E_{il}\t u ) \mA_n=
\mA_1E_{il}\mA_1\t \mA_{n-1} u \mA_{n-1} = F\t \mA_{n-1} =
\gp_1.$$ By symmetry, all  $\gp_i\subseteq \mA_na \mA_n$, and so
$\ga_n = \gp_1+\cdots + \gp_n\subseteq \mA_n a\mA_n$. $\ga_n$ is
 the maximal  ideal of $\mA_n$ that is properly contained in
 the ideal $\mA_na\mA_n$ (since $a$ is generic), and so $\mA_n
 a\mA_n = \mA_n$, as required.

8. Let, for a moment, $n=1$. One can easily verify that, for all
$i,j,k\in \N$, 
\begin{equation}\label{lxkdk}
x^kE_{ij} = E_{i+k, j}, \;\; \der^k E_{ij} = i(i-1) \cdots (
i-k+1) E_{i-k, j},
\end{equation}
\begin{equation}\label{rxkdk}
 E_{ij} x^k = E_{i,j-k}, \;\; E_{ij} \der^k =
(j+k) (j+k-1) \cdots ( j+1) E_{i,j+k},
\end{equation}
where $E_{st}:=0$ if either $s<0$ or $t<0$. By (\ref{lxkdk}), for
each $j\in \N$, the left $\mA_1$-module $C_j:= \oplus_{i\in
\N}KE_{ij}$ is isomorphic to the left $\mA_1$-module $K[x]$, and
so $C_j$ is a faithful, simple, left $\mA_1$-module. The left
$\mA_1$-module 
\begin{equation}\label{lF}
F=\oplus_{i\in \N}\, C_i\simeq K[x]^{(\N )}
\end{equation}
is a direct sum of $\N$ copies of $K[x]$, and so ${}_{\mA_1}F$ is
the {\em faithful, semi-simple, left } $\mA_1$-module.

For an arbitrary $n$, the left $\mA_n$-module $C_{i_1}\t \cdots \t
C_{i_n}$ is isomorphic to $P_n$. Therefore, 
\begin{equation}\label{1lF}
{}_{\mA_n}F^{\t n } \simeq \oplus_{i_1, \ldots, i_n\in \N}\,
 C_{i_1}\t \cdots \t C_{i_n}\simeq P_n^{(\N^n)}
\end{equation}
is a semi-simple, left $\mA_n$-module which is faithful, by
statement 6.

Recall that the structure of the right $\mA_n$-module $P_n$ is
given by the rule: $p*a:= \th (a) p$ where $p\in P_n$ and $a\in
\mA_n$. By (\ref{thFn}), $\th (F^{\t n } ) = F^{\t n }$. Then, by
(\ref{1lF}), the $K$-module $F^{\t n } \simeq P_n^{(\N^n)}$ has
the natural structure of {\em right} $\mA_n$-module, namely,
$f*a:= \th (a) f$ where $f\in F^{\t n }$ and $a\in \mA_n$. In
order to distinguish this structure of right $\mA_n$-module from
the obvious structure (as a right ideal of $\mA_n$) we write
$F^{\t n }_\th \simeq P_n^{(\N^n)}$. The map 
\begin{equation}\label{FthPn}
\th : F^{\t n } \ra F^{\t n } _\th \simeq P_n^{(\N^n)}, \;\;
f\mapsto \th (f),
\end{equation}
is an {\em isomorphism} of right $\mA_n$-modules (since $\th (fa)
= \th (a) \th (f) = \th (f) *a$). Therefore, $F^{\t n }_{\mA_n}
\simeq F^{\t n}_\th \simeq P_n^{(\N^n)}$ is the faithful,
semi-simple, right $\mA_n$-module by the proved left version of
this fact.

The $\mA_n$-bimodule $F^{\t n }$ is simple since the ring $F^{\t n
} \simeq M_\infty (K)^{\t n } \simeq M_\infty (K)$ is simple. The
map $1\t \th : \mA_n \t \mA_n^{op} \ra \mA_{2n}$ is an isomorphism
of $K$-algebras such that $1\t \th (F^{\t n } \t F^{\t n }) =
F^{\t 2n}$. Using this equality and ${}_{\mA_n}F^{\t n
}_{\mA_n}\simeq {}_{\mA_n\t \mA_n^{op}}F^{\t n }\simeq
{}_{\mA_{2n}}F^{\t n }$, we see that the $\mA_n$-bimodule $F^{\t n
}$ is faithful: if $\ga F^{\t n }=0$ for some nonzero ideal $\ga$
of $\mA_{2n}$ then $F^{\t 2n} \subseteq \ga$, and so $0= \ga F^{\t
n } \supseteq F^{\t 2n} \cdot F^{\t n } = F^{\t n } F^{\t n }
F^{\t n } = F^{\t n } \neq 0$, a contradiction.

9. Let $\soc (\mA_n)$ be the socle of the module ${}_{\mA_n}\mA_n$
(resp. ${}_{\mA_n}{\mA_n}_{\mA_n}$). By statement 8, $F^{\t n }
\subseteq \soc (\mA_n)$. Suppose that $F^{\t n } \neq \soc
(\mA_n)$. Then $\soc (\mA_n) = F^{\t n } \oplus M$ for a nonzero
module ${}_{\mA_n}M$ (resp. ${}_{\mA_n}M_{\mA_n}$). On the one
hand, $F^{\t n } \soc (\mA_n) \subseteq F^{\t n } $ and so $F^{\t
n } M=0$. On the other hand, $F^{\t n } M\neq 0$, by statement 6,
a contradiction. Therefore, $\soc (\mA_n) = F^{\t n }$.

The algebra $\mA_n$ admits the involution $\th $ such that $\th
(F^{\t n }) = F^{\t n }$. Therefore, $\soc ({\mA_n}_{\mA_n}) =
F^{\t n }$ since $\soc ({}_{\mA_n} \mA_n) = F^{\t n }$.

10. Let $M$ be a faithful, simple, left (resp. right)
$\mA_n$-module. Then $F^{\t n } M\neq 0$ (resp. $MF^{\t n } \neq
0$). Choose a nonzero element, say $m$, of $M$ such that $F^{\t n
} m\neq 0$ (resp. $mF^{\t n } \neq 0$). Then $M= F^{\t n } m$
(resp. $M= mF^{\t n }$), by simplicity of $M$. There is the
epimorphism $F^{\t n } \ra F^{\t n }m= M$, $f\mapsto fm$ (resp.
$F^{\t n } \ra mF^{\t n }= M$, $f\mapsto mf$) of left (resp.
right) $\mA_n$-modules. Now, the result follows from statement 8.

11. $3n= \GK ( \mA_n) \geq \GK (\mA_n / \ga ) \geq \GK (\mA_n/
\ga_n ) = \GK (\CA_n) = 3n$, and so $\GK (\mA_n / \ga ) = 3n$.
  $\Box $


\section{Unique factorization of ideals of $\mA_n$ and $\Spec (\mA_n)$}\label{ImAnSP}

In this section, all the results on ideals that are mentioned in
the Introduction are proved.

Let $\CB_n$ be the set of all functions $f:\{ 1, 2, \ldots , n\}
\ra \Z_2:= \{ 0,1\}$ where $\Z_2:= \Z / 2\Z$ is a field. $\CB_n$
is a commutative ring with respect to addition and multiplication
of functions.  For $f,g\in \CB_n$, we write $f\geq g$ iff $f(i)
\geq g(i)$ for all $i=1, \ldots , n$ where $1>0$. Then $(\CB_n,
\geq )$ is a partially ordered set. For each function $f\in
\CB_n$, $I_f$ denotes the ideal $I_{f(1)}\t \cdots \t I_{f(n)}$ of
$\mA_n$ which is the tensor product of the ideals $I_{f(i)}$  of
the tensor components $\mA_1(i)$ in $\mA_n= \mA_1(1)\t \cdots \t
\mA_1(n)$ where $I_0:= F$ and $I_1:= \mA_1$. $f\geq g$ iff
$I_f\supseteq I_g$. For $f,g\in \CB_n$, $I_fI_g= I_f\cap I_g =
I_{fg}$. By induction on the number of functions one immediately
proves that, for $f_1, \ldots , f_s\in \CB_n$,
$$\prod_{i=1}^s I_{f_i}= \cap_{i=1}^s I_{f_i} = I_{f_1\cdots
f_s}.$$ Let $\CC_n$ be the set of all  subsets of $\CB_n$ all
distinct elements of which are incomparable (two distinct elements
$f$ and $g$ of $\CB_n$ are {\em incomparable} iff $f\not\leq g$
and $g\not\leq f$). For each $C\in \CC_n$, let $I_C:= \sum_{f\in
C}I_f$, the ideal of $\mA_n$. The next result classifies ideals of
the algebra $\mA_n$.

\begin{theorem}\label{C15Ma7}
Let $K$ be a field of characteristic zero. Then
\begin{enumerate}
\item The map $C\mapsto I_C:= \sum_{f\in C}I_f$ from the set
$\CC_n$ to the set of  ideals of $\mA_n$ is a bijection where
$I_\emptyset :=0$. In particular, there are only finitely many
ideals of $\mA_n$. \item Each ideal $I$ of $\mA_n$ is an
idempotent ideal, i.e. $I^2= I$. \item Ideals of $\mA_n$ commute
($IJ= JI$).
\end{enumerate}
\end{theorem}

{\it Proof}. 1. Statement 1 follows from Lemma \ref{a27Ma7}.

2. The result is obvious for $I=0$. So, let $I\neq 0$. By
statement 1, $I= \sum_{f\in C} I_f$ for some $C\in \CC_n$. Then
$I^2= \sum_{f\in C}I^2_f+\sum_{f\neq g} I_fI_g=\sum_{f\in
C}I_f+\sum_{f\neq g} I_fI_g=\sum_{f\in C}I_f=I$.

3. $I_fI_g = I_{fg} = I_{gf}= I_gI_f$ for all $f,g \in \CB_n$. The
result is obvious if either  $I=0$ or $J=0$. So, let $I\neq 0$ and
$J\neq 0$. By statement 1, $I= I_C$ and $J= I_D$ for some $C,D\in
\CC_n$. Then $IJ= (\sum_{f\in C} I_f)(\sum_{g\in D} I_g )
=(\sum_{g\in D} I_g)(\sum_{f\in C} I_f ) =JI$. $\Box$

Let $B_n$ be the $K$-basis for the algebra $\mA_n$ that is the
tensor product of the $K$-bases from Theorem \ref{12Ma7}.(2) (see
the proof of Corollary \ref{15Ma7}.(7) for details). For each
element $b= b_1\t \cdots \t b_n$ of $B_n$, one can attach an
element $f_b$ of $\CB_n$ by the rule
$$
f_b(i):=\begin{cases}
1& \text{if $b_i\not\in F$},\\
0& \text{if $b_i\in F$}.\\
\end{cases}
$$
Let $a= \sum_{b\in B_n} \l_b b$ be a nonzero element of $\mA_n$
where $\{ \l_b \}$ are the coordinates of $a$ with respect to the
basis $B_n$. One has the well-defined map $\mA_n  \ra \CC_n$,
$a\mapsto \Max (a)$, where $\Max (a)$ are the maximal elements (in
$\CB_n$) of the subset $\{ f_b \, | \, \l_b\neq 0\}$ of $\CB_n$
where $\Max (0):=\emptyset$.

\begin{lemma}\label{a27Ma7}
Let $K$ be a field of characteristic zero and $ 0\neq a=
\sum_{b\in B_n} \l_bb\in \mA_n$. Then $\mA_n a\mA_n = \sum_{b\in
\Max (a) } I_{f_b}= \sum_{\{ b \, | \, \l_b\neq 0 \} } I_{f_b}$.
\end{lemma}

{\it Proof}. It suffices to prove only the first equality since
the second follows from the first: for any $b\in B_n$ such that
$\l_b\neq 0$ there exists $c\in \Max (a)$ such that $\l_c\neq 0$
and $f_c\geq f_b$, and so $I_{f_c}\supseteq I_{f_b}$. Hence,
$\sum_{b\in \Max (a) } I_{f_b}= \sum_{\{ b \, | \, \l_b\neq 0 \} }
I_{f_b}$.

Fix $b\in \Max ( a)$. Up to order of the tensor multiples in the
tensor product  $\mA_n= \t_{i=1}^n \mA_1(i)$, one may assume that
$$
f_b(i):=\begin{cases}
1& \text{if $1\leq i \leq s$},\\
0& \text{if $s+1\leq i \leq n$}.\\
\end{cases}
$$
We have to show that $\mA_n a \mA_n \supseteq \mA_s\t F^{\t
(n-s)}$. $a= \l_1b_1+\cdots + \l_tb_t+\l_{t+1}b_{t+1} +\cdots +
\l_rb_r$ where $b:=b_1$, all $\l_i\neq 0$, $f_{b_1} =\cdots =
f_{b_t}$ and $f_{b_j} \neq f_b$ for $j=t+1, \ldots , r$. By the
choice of $b$, for each $j$ such that $t+1\leq j \leq r$, either
$f_{b_j} <f_b$ or, otherwise, the functions $f_{b_j}$ and $f_b$
are incomparable.  $b_1= c_1\t f_1$ for unique elements $c_1\in
\CM_s$ and $f_1= E_{p_1, q_1} (s+1) \cdots E_{p_{n-s},
q_{n-s}}(n)\in \CN_{n-s}$ where $\CM_s$ and $\CN_{n-s}$ were
defined in the proof of Corollary \ref{15Ma7}.(7). Let
$E_{pp}:=E_{p_1, p_1} (s+1) \cdots E_{p_{n-s}, p_{n-s}}(n)$ and
$E_{qq}:= E_{q_1, q_1} (s+1) \cdots E_{q_{n-s}, q_{n-s}}(n)$. Then
$E_{pp}f_1E_{qq}=f_1$. Note that $E_{pp} \mA_{n-s} E_{qq} =
KE_{pq}= Kf_1$. Changing the element $a$ for the element
$E_{pp}aE_{qq}$ and deleting zero terms of the type $E_{pp}\l_\nu
b_\nu E_{qq}$, one may assume that $a= (\sum_{i=1}^l\mu_ic_i) \t
f_1$ where all $\mu_i\neq 0$ and $\mu_i\in K$; all $c_i\in B_s$;
$f_{c_1}=\cdots = f_{c_k}=1$ for some $k$ such that $1\leq k \leq
l$, i.e. $c_1, \ldots , c_k\in \CM_s$; and $f_{c_j} <f_{c_1}$
where $k+1\leq j \leq l$. The element
$c:=\sum_{i=1}^l\mu_ic_i\not\in \ga_s$, hence $\mA_sc\mA_s=
\mA_s$. Now, $\mA_na\mA_n = \mA_n (c\t f) \mA_n=\mA_s c \mA_s \t
\mA_{n-s} f_1\mA_{n-s} = \mA_s\t F^{\t (n-s)}$, as required. $\Box
$

The next result is a useful criterion of when one ideal contains
another.

\begin{corollary}\label{i29Ma7}
Let $K$ be a field of characteristic zero and $C, C'\in \CC_n$.
 Then $I_C\subseteq I_{C'}$ iff $C\leq C'$ (this means that, for
 each $f\in C$, there exists $f'\in C'$ such that $f\leq f'$).
\end{corollary}

{\it Proof}. This follows at once from Lemma \ref{a27Ma7} and
Theorem \ref{C15Ma7}.(1). $\Box$

\begin{corollary}\label{b27Ma7}
Let $K$ be a field of characteristic zero and $s_n$ be the number
of  ideals of $\mA_n$. Then $2-n+\sum_{i=1}^n2^{n\choose i } \leq
s_n\leq 2^{2^n}$.
\end{corollary}

{\it Proof}. Let $\Sub_n$ be the set of all subsets of $\{ 1,
\ldots, n \}$. $\Sub_n$ is a partially ordered set with respect to
`$\subseteq $'. For each $f\in \CB_n$, the subset $\supp (f) := \{
i \, | \, f(i) =1\}$ of $\{ 1, \ldots , n\}$ is called the {\em
support} of $f$. The map $ \CB_n\ra \Sub_n$, $ f\mapsto \supp
(f)$,  is an isomorphism of posets. Let $\SSub_n$ be the set of
all subsets of $\Sub_n$. An element $\{ X_1, \ldots , X_s\}$ of
$\SSub_n$ is called {\em incomparable} if for all $i\ne j$ such
that $1\leq i,j\leq s$ neither $X_i\subseteq X_j$ nor
$X_i\supseteq X_j$. An empty set and one element set are called
incomparable by definition. Let $\Inc_n$ be the subset of
$\SSub_n$ of all incomparable elements of $\SSub_n$. Then the map
\begin{equation}\label{CnInc}
\CC_n\ra \Inc_n, \;\; \{ f_1, \ldots , f_s\} \mapsto \{ \supp
(f_1), \ldots , \supp (f_s)\},
\end{equation}
is a bijection. For each $i=1, \ldots , n$, there are precisely
${n\choose i}$ subsets of $\{ 1,\ldots , n\}$ that contain exactly
$i$ elements. Any non-empty collection of these is an incomparable
set, hence $s_n\geq 2+\sum_{i=1}^n(2^{n\choose i} -1)=
2-n+\sum_{i=1}^n 2^{n\choose i}$ where  2 `represents' the zero
ideal and the ideal $F^{\t n }$ which corresponds to an empty set.
Clearly, $s_n= |\Inc_n| \leq 2^{2^n}$.  $\Box $

The next corollary classifies all the prime ideals of $\mA_n$.
\begin{corollary}\label{c27Ma7}
Let $K$ be a field of characteristic zero. Then the map
$$\Sub_n\ra \Spec (\mA_n ), \;\;  I\mapsto \gp_I:=\sum_{i\in
I}\gp_i, \;\; \emptyset \mapsto 0,$$ is a bijection, i.e. any
nonzero prime ideal of $\mA_n$ is a unique sum of height 1 primes;
$|\Spec (\mA_n)|=2^n$; the height of $\gp_I$ is $|I|$.
\end{corollary}

{\it Proof}. 0 is the prime ideal since  $\mA_n$  is a prime ring.
Let $I$ be a nonzero subset of $\{ 1, \ldots , n\}$ that contains,
say $s$, elements. Then $\mA_n/\gp_I\simeq \CA_s\t \mA_{n-s}$. The
ring $\CA_s$ is a central, simple $K$-algebra, hence the map $\ga
\mapsto \CA_s\t \ga$ is a bijection from the set of ideals of
$\mA_{n-s}$ to the set of ideals of  $\CA_s\t \mA_{n-s}$ (this is
an easy consequence of the Density Theorem). Then, $\CA_s\t F^{\t
(n-s)}$ is the smallest nonzero ideal of $\CA_s\t \mA_{n-s}$ and
it is idempotent, and so $\CA_s\t \mA_{n-s}$ is a prime ring. This
means that $\gp_I$ is a prime ideal. By Theorem \ref{C15Ma7}.(1),
the map $\Sub_n\ra \Spec ( \mA_n)$, $I\mapsto \gp_I$, is an
injection. It remains to prove that it is a surjection. Let $\gp$
be a prime nonzero ideal of $\mA_n$. Then $\gp = \sum_{f\in C}I_f$
for some $C\in \CC_n$ (Theorem \ref{C15Ma7}.(1)). If $|\supp (f)|=
n-1$ for all $f\in C$ then $I_f= \gp_i$ where $i=i(f)$ is a unique
element of the set $\{ 1, \ldots , n\}$ such that $f(i) =0$, and
so $\gp = \sum \gp_i$.

Suppose that $|\supp (f)| \neq n-1$ for some $f\in C$. Let us show
that this case is not possible since then $\gp$ would {\em not} be
a prime ideal. One can choose two distinct functions, say $g,h\in
\CB_n$,  such that $g>f$, $h>f$, and $gh= f$. Then $I_f = I_{gh} =
I_gI_h$. Let $\ga:= I_g+\gc$ and $\gb := I_h+\gc$ where $\gc :=
\sum_{f\neq f'\in C} I_{f'}$. The ideals $\ga$ and $\gb$ {\em
strictly} contain the ideal $\gp$ and $$\ga \gb = (I_g+\gc ) (
I_h+\gc ) = I_gI_h+ \gc I_h+ I_g\gc +\gc^2= I_{gh}+ \gc I_h+
I_g\gc +\gc^2\subseteq I_f+\gc = \gp .$$ This contradicts to the
fact that $\gp $ is a prime ideal, and we are done.

It is obvious that $|\Spec (\mA_n)|=2^n$. The  fact that the
height of the ideal $\gp_I$ is $|I|$ follows from Lemma
\ref{p28Ma7}.(1).  $\Box $

The next criterion of when a prime ideal contains another prime is
used in finding the classical Krull dimension of the algebra
$\mA_n$ (Corollary \ref{k28Ma7}).
\begin{lemma}\label{p28Ma7}
Let $K$ be a field of characteristic zero; $\gp , \gq \in \Spec
(\mA_n)$; $ \gp = \gp_{i_1}+\cdots + \gp_{i_s}$ and $\gq =
\gp_{j_1}+\cdots + \gp_{j_t}$ be their decompositions as in
Corollary \ref{c27Ma7}. Then
\begin{enumerate}
\item $\gp \subseteq \gq$ iff $\{ \gp_{i_1}, \ldots , \gp_{i_s}\}
\subseteq \{ \gp_{j_1}, \ldots ,  \gp_{j_t}\} $. \item If $\gp
\subseteq \gq$ then $\gp\gq= \gp$. \item The poset $(\Spec
(\mA_n), \subseteq )$ is an isomorphic to the set $\Sub_n$ of all
subsets of $\{ 1, \ldots n \}$.
\end{enumerate}
\end{lemma}

{\it Proof}. 1. $(\Rightarrow )$ $(\gq + \gp ) / \gp$ is the ideal
of the algebra $\mA_n / \gp \simeq \CA_s\t \mA_{n-s}$. The algebra
$\CA_s$ is central and simple. By the Density Theorem, each ideal
of the algebra $\CA_s\t \mA_{n-s}$ is of the type $\CA_s\t \ga$
for some ideal $\ga$ of  $\mA_{n-s}$. By Corollary \ref{c27Ma7}
(applied to the algebra $\mA_{n-s}$), we have the inclusion $\{
\gp_{i_1}, \ldots , \gp_{i_s}\} \subseteq \{ \gp_{j_1}, \ldots ,
\gp_{j_t}\} $.

$(\Leftarrow )$ If $\{ \gp_{i_1}, \ldots , \gp_{i_s}\} \subseteq
\{ \gp_{j_1}, \ldots ,  \gp_{j_t}\} $ then $ \gp =
\gp_{i_1}+\cdots + \gp_{i_s} \subseteq  \gp_{j_1}+\cdots +
\gp_{j_t}=\gq $.

2. By statement 1, if $\gp \subseteq \gq$ then $\gq = \gp + \gr$
for an ideal $\gr$. Then $\gp \gq = \gp^2+ \gp\gr= \gp + \gp \gr=
\gp$.

3. Statement 3 follows from Corollary \ref{c27Ma7} and statement
1. $\Box $

For each $s=0, 1, \ldots , n$, there are precisely ${n\choose s}$
prime  ideals of height $s$, namely, $$\{ \gp_{i_1}+\cdots
+\gp_{i_s}\, | \, 1\leq i_1<\cdots <i_s\leq n\}.$$

\begin{corollary}\label{k28Ma7}
Let $K$ be a field of characteristic zero. Then the classical
Krull dimension of $\mA_n$ is $n$.
\end{corollary}

{\it Proof}. By Lemma \ref{p28Ma7}.(1), $0\subset \gp_1\subset
\gp_1+\gp_2\subset\cdots \subset \gp_1+\cdots +\gp_n$ is a longest
chain of primes, and so ${\rm cl.K.dim} (\mA_n) = n$. $\Box $

$(\Spec (\mA_n) , \subseteq )$ is a poset. Two primes $\gp $ and
$\gq$ are called {\em incomparable} if neither $\gp \subseteq \gq
$ nor $\gp \supseteq \gq$.

For each ideal $\ga$ of $\mA_n$ such that $\ga\neq \mA_n$, let
$\Min (\ga)$ be the set of all minimal primes over $\ga$. The set
$\Min (\ga ) $ is a non-empty set since the ring $\mA_n$ has only
finitely many primes.

For each $f\in \CB_n$, the set $\csupp (f):= \{ i\, | \, f(i)=0\}$
is called the {\em co-support} of $f$. Clearly, $\csupp (f) = \{
1, \ldots , n\} \backslash \supp (f)$.

\begin{theorem}\label{27Ma7}
Let $K$ be a field of characteristic zero. Then
\begin{enumerate}
\item Each ideal $\ga$ of $\mA_n$ such that $\ga \neq \mA_n$ is a
unique product of incomparable primes, i.e. if $\ga = \gq_1\cdots
\gq_s= \gr_1\cdots \gr_t$ are two such products then $s=t$ and
$\gq_1= \gr_{\s (1)}, \ldots , \gq_s= \gr_{\s (s)}$ for a
permutation $\s$ of $\{ 1, \ldots, n\}$. \item  Each ideal $\ga$
of $\mA_n$ such that $\ga \neq \mA_n$ is a unique intersection of
incomparable primes, i.e. if $\ga = \gq_1\cap \cdots\cap  \gq_s=
\gr_1\cap \cdots \cap \gr_t$ are two such intersections then $s=t$
and $\gq_1= \gr_{\s (1)}, \ldots , \gq_s= \gr_{\s (s)}$ for a
permutation $\s$ of $\{ 1, \ldots, n\}$. \item   For each ideal
$\ga$ of $\mA_n$ such that $\ga \neq \mA_n$, the sets of
incomparable primes in statements 1 and 2 are the same,
$\ga=\gq_1\cdots \gq_s=\gq_1\cap \cdots\cap  \gq_s$. \item The
ideals $\gq_1,\ldots , \gq_s$ in statement 3 are the minimal
primes of $\ga$, and so  $\ga = \prod_{\gp \in \Min (\ga )}\gp
=\cap_{\gp \in \Min (\ga )}\gp$.
\end{enumerate}
\end{theorem}

{\it Proof}. 1. For each ideal $\ga$ of $\mA_n$, we have to prove
that $\ga$ is a product of incomparable primes and that this
product is unique. Since the ring $\mA_n$ is prime these two
statements are obvious when $\ga =0$. So, let $\ga \neq 0$.

{\em Existence}: Let $f\in \CB_n$; then $I_f=\prod_{i\in \csupp
(f)}\gp_i$. Let $\gb$ be any ideal of $\mA_n$. Since $\gb^2= \gb$,
it follows at once that 
\begin{equation}\label{Ifb}
I_f+\gb = \prod_{i\in \csupp (f)} (\gp_i+\gb ).
\end{equation}
By Theorem \ref{C15Ma7}.(1), $\ga = I_{f_1} +\cdots + I_{f_s}$ for
some $f_i\in \CB_n$. Repeating $s$ times (\ref{Ifb}), we see that
\begin{equation}\label{aiSup}
\ga = \prod_{i_1\in \csupp (f_1),\ldots , i_s\in \csupp
(f_s)}(\gp_{i_1} +\cdots + \gp_{i_s})
\end{equation}
is the product of primes, by Corollary \ref{c27Ma7}. Note that
ideals of $\mA_n$ commute; each ideal is an idempotent ideal; and
if $\gp \subseteq \gq$ is an inclusion of primes then $\gp \gq =
\gp$. Using these three facts and (\ref{aiSup}), we see that $\ga$
is a product of {\em incomparable} primes.

Uniqueness follows from the next lemma which will be used several
times in the proof of this theorem.

\begin{lemma}\label{u28Ma7}
Let $\{ \gq_1, \ldots , \gq_s\}$ and $\{ \gr_1, \ldots , \gr_t\}$
be two sets of incomparable ideals of a ring such that each ideal
from the first set contains an ideal from the second and each
ideal from the second set contains an ideal from the first. Then
$s=t$ and $\gq_1= \gr_{\s (1)}, \ldots , \gq_s= \gr_{\s (s)}$ for
a permutation $\s$ of $\{ 1, \ldots, n\}$.
\end{lemma}

{\it Proof Lemma \ref{u28Ma7}}. For each $\gq_i$, there are ideals
$\gr_j$ and $\gr_k$ such that $\gr_j\subseteq \gq_i\subseteq
\gr_k$, hence $\gq_i= \gr_j = \gr_k$ since the ideals $\gr_j$ and
$\gr_k$ are incomparable if distinct. This proves that for each
ideal $\gq_i$ there exists a unique ideal, say $\gr_{\s (i)}$,
such that $\gq_i= \gr_{\s (i)}$. By symmetry, for each ideal
$\gr_j$ there exists a unique  ideal, say $\gq_{\tau (j)}$, such
that $\gr_j = \gq_{\tau (j)}$. Then, $s=t$ and $\gq_1= \gr_{\s
(1)}, \ldots , \gq_s= \gr_{\s (s)}$ for the  permutation $\s$ of
$\{ 1, \ldots, n\}$. $\Box $

{\em Uniqueness}: Let $\ga = \gq_1\cdots \gq_s= \gr_1\cdots \gr_t$
 be  two products of incomparable primes. Each ideal $\gq_i$ contains
an   ideal $\gr_j$, and each ideal $\gr_k$ contains an ideal
$\gq_l$. By Lemma \ref{u28Ma7}, $s=t$ and $\gq_1= \gr_{\s (1)},
\ldots , \gq_s= \gr_{\s (s)}$ for a permutation $\s$ of $\{ 1,
\ldots, n\}$.

2. {\em Uniqueness}: Suppose that an ideal $\ga$ has two
presentations $\ga = \gq_1\cap \cdots \cap\gq_s= \gr_1\cap
\cdots\cap  \gr_t$ of incomparable primes. The sets $\{ \gq_1,
\ldots , \gq_s\}$ and $\{ \gr_1, \ldots , \gr_t\}$ of incomparable
primes satisfy the conditions of Lemma \ref{u28Ma7}, and so
uniqueness follows.

{\em Existence}: Let $\CI$ be the set of all the ideals of
$\mA_n$, and $\CI'$ be the set of ideals of $\mA_n$ that are
intersection of incomparable primes. Then $\CI'\subseteq \CI$. The
map
$$\CI \ra \CI',\;\; \gq_1\cdots \gq_s\mapsto \gq_1\cap \cdots \cap
\gq_s,$$ is a  bijection since $|\CI |<\infty$ and by uniqueness
of presentations $\gq_1\cdots \gq_s$ (statement 1) and $\gq_1\cap
\cdots \cap \gq_s$ (see above) where $\gq_1, \ldots , \gq_s$ are
incomparable primes. Then $\CI = \CI'$. This proves that each
ideal $\ga$ of $\mA_n$ is  an intersection of incomparable primes.

3. Let $\ga$ be an ideal of $\mA_n$ and $\ga = \gq_1\cdots \gq_s=
\gr_1\cap \cdots \cap \gr_t$ where $S:=\{ \gq_1, \ldots , \gq_s\}$
and $T:=\{ \gr_1, \ldots , \gr_t\}$ are  sets of incomparable
primes. The sets $S$ and $T$ satisfy the conditions of Lemma
\ref{u28Ma7}, and so $s=t$ and $\gq_1= \gr_{\s (1)}, \ldots ,
\gq_s= \gr_{\s (s)}$ for a permutation $\s$ of $\{ 1, \ldots,
n\}$. This means that $\ga = \gq_1\cdots \gq_s= \gq_1\cap \cdots
\cap \gq_s$.

4. Let $\ga = \gq_1\cdots \gq_s = \gq_1\cap \cdots \cap \gq_s$ be
as in statement 3 and let $\Min (\ga )=\{ \gr_1, \ldots , \gr_t\}$
be the set of minimal primes over $\ga$. Then $\Min (\ga )
\subseteq S:=\{ \gq_1, \ldots \gq_s\}$ ($\ga = \gq_1\cdots
\gq_s\subseteq \gr_i$ implies $\gq_j\subseteq \gr_i$ for some $j$,
and so $\gq_j = \gr_i$ by the minimality of $\gr_i$). Up to order,
let $\gr_1= \gq_1, \ldots , \gr_t = \gq_t$. It remains to show
that $t=s$. Suppose that $t<s$, we seek a contradiction. This
means that each prime $\gq_i$, $i=t+1, \ldots , s$, contains $\ga$
and is {\em not} minimal over $\ga$. Hence, $\gq_i$ contains a
minimal prime, say $\gq_{\tau (i)}$, a  contradiction (the ideal
$\gq_i$ and $\gq_{\tau (i)}$ are incomparable).  $\Box $

\begin{corollary}\label{a28Ma7}
Let $K$ be a field of characteristic zero, $\ga$ and $\gb$ be
ideals of $\mA_n$ distinct from $\mA_n$ in statement 1, 2 and 5.
Then
\begin{enumerate}
\item $\ga = \gb$ iff $\Min (\ga ) = \Min (\gb )$.\item $\Min (\ga
\cap \gb ) = \Min (\ga \gb ) = $ the set of minimal elements (with
respect to inclusion) of the set $\Min (\ga ) \cup \Min (\gb )$.
\item $ \ga\cap \gb = \ga \gb$. \item If $\ga \subseteq \gb$ then
$\ga \gb = \ga$. \item $\ga \subseteq \gb$ iff $\Min (\ga )
 \eqslantless \Min (\gb )$ (the $\eqslantless$ means that
and each $\gq\in \Min (\gb )$ contains some $\gp\in \Min (\ga )$).
\end{enumerate}
\end{corollary}

{\it Proof}. 1. Statement 1 is obvious due to Theorem
\ref{27Ma7}.(4).

2. Let $\CM$ be the set of minimal elements of the union $\Min
(\ga ) \cup \Min (\gb )$. The elements of $\CM$ are incomparable,
and (by Theorem \ref{27Ma7}.(4))
$$ \ga \cap \gb = \cap_{\gp\in \Min (\ga )}\cap \cap_{\gq \in \Min
(\gb )}\gq = \cap_{\gr \in \CM } \gr.$$ By Theorem
\ref{27Ma7}.(2), $\Min (\ga \cap \gb ) = \CM$. By Lemma
\ref{p28Ma7}.(2), $$\ga \gb = \prod_{\gp \in \Min (\ga ) }\gp
\cdot \prod_{\gq \in \Min (\gb )}\gq = \prod_{\gr \in \CM}\gr =
\ga \cap \gb.$$

3. The result is obvious if one the ideals is equal to $\mA_n$.
So, let the ideals are distinct from $\mA_n$.  By statement 2,
$\Min (\ga \cap \gb ) = \Min (\ga \gb )$, then, by statement 1,
$\ga\cap \gb = \ga \gb$.

4. If $\ga \subseteq \gb$ then, by statement 3, $\ga \gb = \ga\cap
\gb = \ga$.

5. $(\Rightarrow )$ If $\ga \subseteq \gb$ then $\Min (\ga )
\eqslantless\Min (\gb )$ since $\ga = \prod_{\gp\in \Min (\ga
)}\gp \subseteq \prod_{\gq \in \Min (\gb )}\gq = \gb$.

$(\Leftarrow )$ Suppose that $\Min (\ga ) \eqslantless\Min (\gb
)$. For each $\gq\in \Min (\gb )$, let $S(\gq )$ be the set
(necessarily nonempty) of $\gp\in \Min (\ga )$ such that $\gp
\subseteq \gq$. Then $\Min (\ga ) \supseteq S:= \cup_{\gq\in \Min
(\gb )} S(\gq )$ and
$$\ga =\cap_{\gp \in \Min (\ga )}\gp \subseteq \cap_{\gp \in S}\gp \subseteq
\cap_{\gq\in \Min (\ga )}\gq = \gb.\;\;\; \Box $$

\begin{theorem}\label{28Ma7}
Let $K$ be a field of characteristic zero. Then the lattice of
ideals of the algebra $\mA_n$ is distributive, i.e. $(\ga \cap \gb
) \gc = \ga \gc\cap \gb \gc$ for all ideals $\ga$, $\gb$, and
$\gc$.
\end{theorem}

{\it Proof}.  By Corollary \ref{a28Ma7}.(3), $(\ga \cap \gb ) \gc
= \ga \cap \gb \cap \gc = (\ga \cap \gc ) \cap ( \gb \cap \gc ) =
\ga \gc \cap \gb \gc$.  $\Box $

\begin{theorem}\label{1Jun7}
Let $K$ be a field of characteristic zero, $\ga$ be an ideal of
$\mA_n$, and $\CM$ be the minimal elements with respect to
inclusion of a set of ideals $\ga_1, \ldots , \ga_k$ of $\mA_n$.
Then
\begin{enumerate}
\item  $\ga = \ga_1\cdots \ga_k$ iff $\Min (\ga ) = \CM$.\item
$\ga = \ga_1\cap \cdots\cap \ga_k$ iff $\Min (\ga ) = \CM$.
\end{enumerate}
\end{theorem}

{\it Proof}. By Corollary \ref{a28Ma7}.(3), it suffices to prove,
say, the first statement.

$(\Rightarrow )$ Suppose that $\ga = \ga_1\cdots \ga_k$ then, by
Theorem \ref{27Ma7}.(4) and Corollary \ref{a28Ma7}.(4), $$\ga =
\prod_{i=1}^k\prod_{\gq_{ij}\in \Min (\ga_i)} \gq_{ij} =
\prod_{\gq \in \CM} \gq,$$
 and so $\Min (\ga ) = \CM$, by Theorem
\ref{27Ma7}.(4).

$(\Leftarrow )$ If $\Min (\ga ) = \CM$ then, by Corollary
\ref{a28Ma7}.(4), $\ga = \ga_1\cdots \ga_k$.  $\Box $

{\bf The involution $c$}. Let $K$ be a field of characteristic
zero  and $\CI (\mA_n)$ be the set of all ideals of the algebra
$\mA_n$. Consider the map $\CC_n\backslash \{\emptyset \} \ra
\CC_n\backslash \{ \emptyset \}$, $C\mapsto C+1$, where for $C=\{
f_1, \ldots , f_s\}$, $C+1:= \{ f_1+1, \ldots , f_s+1\}$. The map
is well-defined: $C\in \CC_n$ iff $\{ \supp (f_1) , \ldots , \supp
(f_s)\} \in \Inc_n$ iff $\{ \csupp (f_1) , \ldots , \csupp (f_s)\}
\in \Inc_n$ where $\csupp (f_i):= \{ 1, \ldots , n\} \backslash
\supp (f_i)$ iff $C+1\in \Inc_n$. Consider the map
$$ c: \CI (\mA_n) \ra \CI (\mA_n), \;\; I_C\mapsto I_{C+1}, \;\;
C\in \CC_n,$$ where $c(0):= 0$. Then, for $C\in \CC_n$,
$c(\sum_{f\in C}I_f)= \sum_{f\in C} c(I_f)$. Note that $c(\mA_n) =
F^{\t n }$, $c(\gp_i) = \prod_{j\neq i } \gp_j$, and $c(\ga_n ) =
\sum_{i=1}^n c(\gp_i)$.

Let $C, C'\in \CC_n$, we write $C\preceq C'$ if for each $f\in C$
there exists $f'\in C'$ such that $f\leq f'$, and for each $g'\in
C'$ there exists $g\in C$ such that $g\leq g'$.

\begin{lemma}\label{c29Ma7}
 Let $K$ be a field of characteristic
zero. Then
\begin{enumerate}
\item $c: \CI (\mA_n) \ra \CI (\mA_n)$ is an involution ($c^2=
{\rm id}$) such that $f\leq g$ implies $c(I_f) \supseteq c(I_g)$.
\item $c( \ga ) =\ga$ iff $\ga =I_C$ for some $C=\{ f_1, f_1+1,
\ldots , f_s, f_s+1\}$. \item If $C,C'\in \CC_n$ and $C\preceq C'$
then $c(I_C)\supseteq c(I_{C'})$. \item $c (\ga + \gb ) \subseteq
c(\ga ) +  c(\gb )$ for all ideals $\ga$ and $\gb $ of $\mA_n$. If
$\ga = I_C$ and $\gb = I_{C'}$ for some $C, C'\in \CC_n$ such that
$C\cup C'\in \CC_n$ then $c(\ga + \gb ) = c(\ga ) + c(\gb )$.
\item $c( \ga \gb ) \supseteq c(\ga ) + c(\gb )$ for all nonzero
ideals $\ga $ and $\gb $ of $\mA_n$.
\end{enumerate}
\end{lemma}

{\it Proof}. 1. $c^2= {\rm id}$ since, for all $\emptyset \neq
C\in \CC_n$, $C+1+1= C$. The rest is obvious.

2. $c(I_C ) =I_C$ iff $C+1= C$ iff $C=\{ f_1, f_1+1, \ldots , f_s,
f_s+1\}$.

3. If $C\preceq C'$ then $C'+1\preceq C+1$, and so
$c(I_C)\supseteq c(I_{C'})$.

4. The second statement is obvious: If $C\cup C'\in \CC_n$ then
$$ c(\ga + \gb ) =c(I_{C\cup C'})= \sum_{f\in C\cup C'} I_{f+1}=
\sum_{f\in C} I_{f+1} +\sum_{f\in C'} I_{f'+1}= \ga +\gb .$$ For
arbitrary $\ga$ and $\gb$, let $\ga = I_C$ and $\gb = I_D$, and
$\ga +\gb = I_E$ for some $C,D,E\in \CC_n$. Then, for each $e\in
E$, either $e\geq f$ for some $f\in C$ or $e\geq g$ for some $g\in
D$. Then, either $c(I_e) \subseteq c(I_f)$ or $c(I_e)\subseteq
c(I_g)$. Hence, $c(\ga + \gb ) \subseteq c(\ga ) + c(\gb )$.

5. Let $\ga = I_C$, $\gb = I_D$, and $\ga\gb = I_E$ for some
$C,D,E\in \CC_n$. Then $E\preceq C$ and $E\preceq D$. By statement
3, $c(\ga \gb ) \supseteq c(\ga ) +c(\gb )$. $\Box $

{\bf The involution $\tau$}. Let $K$ be a field of characteristic
zero. For each $\gp \in \Spec (\mA_n)$, there exists a unique
prime ideal $\tau (\gp )$ such that $\ga_n = \gp \oplus \tau (\gp
)$. In more detail, if $\gp = \gp_{i_1} +\cdots + \gp_{i_s}$ then
$\tau (\gp ) =\gp_{j_1} +\cdots + \gp_{j_t}$ where $\{ j_1, \ldots
, j_t\}:= \{ 1, \ldots , n\} \backslash \{ i_1, \ldots , i_s\}$.
The map $\tau : \Spec (\mA_n) \ra \Spec (\mA_n)$ is an order
reversion involution, i.e. $\gp \subseteq \gq$ implies $\tau (\gp
) \supseteq \tau (\gq )$ for $\gp , \gq \in \Spec (\mA_n)$; and
$\tau^2= {\rm id}$. In particular, $\tau $ is an anti-automorphism
of the poset $\Spec (\mA_n)$. $\tau (\ga_n) = 0$ and $\tau (\gp_i)
= \sum_{j\neq i} \gp_j$.  Let $\CI_n$ be the set of ideals of
$\mA_n$ distinct from $\mA_n$. The map $\tau $ can be extended to
the map
$$ \tau : \CI_n\ra \CI_n,\;\; \ga = \cap_{\gq\in \Min (\ga
)}\mapsto \tau (\ga ) := \cap_{\gq\in \Min ( \ga )} \tau (\gq ).$$

\begin{lemma}\label{c28Ma7}
Let $K$ be a field of characteristic zero, and $\ga , \gb \in
\CI_n$. Then
\begin{enumerate}
\item $\tau : \CI_n\ra \CI_n$ is the involution
 ($\tau^2= {\rm id}$). \item  $\tau (\ga ) =\ga$ iff $\tau (\Min
 (\ga ))= \Min (\ga )$.
\end{enumerate}
\end{lemma}

{\it Proof}. 1. The elements of the set $\Min (\ga )$ are
incomparable, then so are the elements of the set $\tau (\Min (\ga
))$, hence $\tau ( \Min (\ga ))= \Min (\tau (\ga ))$. Then
$\tau^2= {\rm id}$.

2. By Corollary \ref{a28Ma7}.(1), $\tau (\ga ) =\ga$ iff $\tau
(\Min (\ga ))= \Min (\ga )$.  $\Box $

A prime ideal of a ring $R$ is called a {\em completely prime} if
$R/\gp$ is a domain.

\begin{corollary}\label{b28Ma7}
Let $K$ be a field of characteristic zero. Then
\begin{enumerate}
\item  $\ga_n$ is the only completely prime ideal of $\mA_n$.
\item $\ga_n$ is the only ideal $\ga$ of $\mA_n$ such that the
factor ring $\mA_n/\ga$ is Noetherian (resp. left Noetherian,
resp. right Noetherian).
\end{enumerate}
\end{corollary}

{\it Proof}. 1. By Corollary \ref{c27Ma7}, any prime ideal $\gp$
of $\mA_n$  is a unique sum $\gp = \gp_{i_1} +\cdots + \gp_{i_s}$.
Then $\mA_n/ \gp\simeq \CA_s\t \mA_{n-s}$. The ring $\CA_s\t
\mA_{n-s}$ is a domain iff $s=n$, that is  $\gp = \ga_n$.

2. The factor ring $\mA_n / \ga_n\simeq \CA_n$ is Noetherian. It
remains to show that $\mA_n/ \ga$ is not left and right Noetherian
for all ideals $\ga$ distinct from $\ga_n$. By Theorem
\ref{27Ma7}.(4), if $\ga \neq \ga_n$ then $\ga_n\not\in \Min (\ga
)$. Choose $\gp \in \Min (\ga )$. Then $\mA_n/\gp \simeq \CA_s\t
\mA_{n-s}$ for some $s\geq 1$. The ring $\CA_s\t \mA_{n-s}$ is not
left or right Noetherian. The ring $\CA_s\t \mA_{n-s}$ is a factor
ring of $\mA_n/\ga$, and the result follows.
 $\Box $


\section{ The group of units $\mA_n^*$ of $\mA_n$ }\label{}

In this section, $K$ be a field of characteristic zero. Let
$\mA_n^*$, $\CA_n^*$ and $K^*$ be the groups of units of the
algebras $\mA_n$, $\CA_n$ and $K$ respectively. Using the
$\Z^n$-grading of the skew polynomial algebra $\CA_n$ (see
(\ref{Anskewlaurent})), it follows that 
\begin{equation}\label{stCAn}
\CA_n^*= (S_n^{-1}\CP_n)^*=\{ K^*\prod_{i\in
\Z}\prod_{j=1}^n(H_j+i)^{n_{ij}}\, | \, (n_{ij})\in (\Z^n)^{(\Z
)}\} \simeq K^*\times (\Z^n)^{(\Z )}
\end{equation}
where the abelian group $(\Z^n)^{(\Z )}$ is the direct sum of $\Z$
copies of $\Z^n$.

{\bf The group $K^*\times \CH_n$}. Let, for a moment, $n=1$. In
this case we usually drop the subscript 1. For each integer $i\geq
1$ and $\l \in K^*$, the element $(H-i)_\l : = H-i+\l\pi_{i-1}$ is
a unit of the algebra $\mD_1$ and its inverse is equal to
$$ (H-i)_\l^{-1}:=
\begin{cases}
\rho_{11}+\l^{-1}\pi_0& \text{if $i=1$},\\
\rho_{1i}+\sum_{j=0}^{i-2} \frac{1}{j+1-i}\, \pi_j +\l^{-1}\pi_{i-1}& \text{if $i\geq 2$}.\\
\end{cases}
$$
As a function of the discrete argument $H$,   $(H-i)_\l^{-1}$
coincides with $\frac{1}{H-i}$ but  instead of having pole at
$H=i$ it takes the value $\l^{-1}$. Consider the following
subgroup of $\mD_1^*$, 
\begin{equation}\label{defH1}
\CH := \{ \prod_{i\geq 0} (H+i)^{n_i}\cdot \prod_{i\geq
1}(H-i)^{n_{-i}}_1\, | \, (n_i)\in \Z^{(\Z )}\}\simeq \Z^{(\Z )}.
\end{equation}
For an arbitrary $n\geq 1$, recall that
$\mA_n=\t_{i=1}^n\mA_1(i)=\mA_1^{\t n }$. For each tensor multiple
$\mA_1(i)= \mA_1$, let $\CH (i)$ be the corresponding group $\CH$.
Their product $\CH_n:= \CH (1) \cdots \CH (n)$ is a subgroup of
$\mD_n^*$ and $\CH_n\simeq  \CH^n\simeq (\Z^n)^{(\Z )}$. The
natural inclusion $\CH_n\simeq (\CH_n+\ga_n)/\ga_n\subset \mA_n/
\ga_n \simeq \CA_n$ and (\ref{stCAn}) yield the isomorphism of
groups 
\begin{equation}\label{KHiA}
K^*\times \CH_n\ra \CA_n^*, \;\; \l\mapsto \l, \;\; H_s+i\mapsto
H_s+i, \;\; (H_s-j)_1\mapsto H_s-j,
\end{equation}
where $\l \in K^*$, $1\leq s\leq n$,  $i\in \N$ and $1\leq j\in
\N$. $K^*\CH_n$ is the subgroup of $\mD_n^*$ such that
$K^*\CH_n\simeq K^*\times \CH_n$.

{\bf The group $(1+F^{\t n })^*$ of units of the monoid $1+F^{\t
n} $}. We are going to find the group $(1+F^{\t n })^*$ of units
of the multiplicative (noncommutative) monoid $1+F^{\t n }$. Let,
for a moment, $n=1$. The ring $F=\oplus_{i,j\in \N} KE_{ij}$ is
the union $M_\infty (K) := \cup_{d\geq 1}M_d(K)= \varinjlim
M_d(K)$ of the matrix algebras $M_d(K):= \oplus_{1\leq i,j\leq
d-1}KE_{ij}$, i.e. $F= M_\infty (K)$.

For each $d\geq 1$, consider the (usual) determinant $\det_d=\det
: 1+M_d(K)\ra K$, $u\mapsto \det (u)$. These determinants
determine the (global) determinant
$$ \det : 1+M_\infty (K)= 1+F\ra K, \;\; u\mapsto \det (u),
$$
where $\det (u)$ is the common value of all determinants
$\det_d(u)$, $d\gg 1$. The (global) determinant has usual
properties of the determinant. In particular, for all $u,v\in
1+M_\infty (K)$, $\det (uv) = \det (u) \cdot \det (v)$. It follows
from  Cramer's formula that the group $\GL_\infty (K):=
(1+M_\infty (K))^*$ of units of the monoid $1+M_\infty (K)$ is
equal to 
\begin{equation}\label{GLiK}
GL_\infty (K) = \{ u\in 1+M_\infty (K) \, | \, \det (u) \neq 0\}.
\end{equation}
Therefore, 
\begin{equation}\label{1GLiK}
(1+F)^* = \{ u\in 1+F \, | \, \det (u) \neq 0\}=GL_\infty (K).
\end{equation}
 The kernel
$$\SL_\infty (K):= \{ u\in
\GL_\infty (K)\, | \, \det (u) =1\}$$
 of the group epimorphism
$\det : \GL_\infty (K)\ra K^*$ is a {\em normal} subgroup of
$\GL_\infty (K)$.

 For any   $n\geq 1$,
\begin{eqnarray*}
F^{\t n } & = & \t_{i=1}^n F(i) = \t_{i=1}^n (\cup_{d_i\geq 1}
M_{d_i}(K))= \cup_{d_1, \ldots , d_n\geq 1} \t_{i=1}^n
M_{d_i}(K)\\
&=& \cup_{d_1, \ldots , d_n\geq 1}M_{d_1 \cdots d_n}(K)= M_\infty
(K).
\end{eqnarray*}
Consider the determinant
$$ \det : 1+F^{\t n} =1+M_\infty (K)\ra K, \;\; u\mapsto \det (u),
$$
as in the case $n=1$. Hence, 
\begin{equation}\label{1pFtn}
(1+F^{\t n })^* = \{ u\in 1+F^{\t n }\, | \, \det (u) \neq 0\} =
(1+M_\infty (K))^*= \GL_\infty (K).
\end{equation}
For each element $u\in (1+F^{\t n })^*$, using Cramer's formula
one can easily find a formula for the inverse $u^{-1}$, it is
Cramer's formula.

\begin{lemma}\label{f21Ma7}
Let $K$ be a field of characteristic zero and $u\in 1+F^{\t n}$.
The following statements are equivalent.
\begin{enumerate}
\item $u\in \mA_n^*$. \item The element $u$ has left inverse in
$\mA_n$ ($vu=1$ for some $v\in \mA_n$). \item The element $u$ has
right  inverse in $\mA_n$ ($uv=1$ for some $v\in \mA_n$). \item
$\det (u)\neq 0$.
\end{enumerate}
\end{lemma}

{\it Proof}. Using the $\Z^n$-grading on $\mA_n$, it is obvious
that the first three statements are equivalent to the fourth.
$\Box $

Since $F^{\t n }$ is an ideal of $\mA_n$, the subgroup $(1+F^{\t n
})^*$ of $\mA_n^*$ is a {\em normal } subgroup: For all $a\in
\mA_n^*$, $a(1+F^{\t n })a^{-1}=1+aF^{\t n }a^{-1} \subseteq
1+F^{\t n }$, and so $a(1+F^{\t n })^*a^{-1} \subseteq (1+F^{\t n
})^*$.

{\bf The subgroup $K^*\times ( \CH_n\ltimes (1+F^{\t n })^*)$ of
$\mA_n^*$}. Let $\mA_n'$ be the subgroup of the group $\mA_n^*$
generated by its subgroups $K^*$, $\CH_n$ and $(1+F^{\t n })^*$.
Let us prove that 
\begin{equation}\label{ApnK}
\mA_n'=K^*\times ( \CH_n\ltimes (1+F^{\t n })^*).
\end{equation}
The subgroup     $(1+F^{\t n })^*$ of $\mA_n^*$ (and of $\mA_n'$)
is normal and the subgroup $K^*$ belongs to the centre of
$\mA_n^*$, hence $\mA_n'= K^*\CH_n(1+F^{\t n })^*$, i.e. each
element $a$ of $\mA_n'$ is a product $a= \l \alpha u$ for some
elements $\l \in K^*$, $\alpha \in \CH_n$ and $u\in (1+F^{\t n
})^*$. In order to prove (\ref{ApnK}) it suffices to show
uniqueness of the decomposition $a= \l \alpha u$. Since $a+\ga_n =
\l \alpha +\ga_n \in (\mA_n/ \ga_n)^* \simeq \CA_n^*$, the
uniqueness of $\l$ and $\alpha$ follows from (\ref{stCAn}) and
(\ref{KHiA}). Then $u= (\l \alpha )^{-1}a$ is unique as well. This
finishes the proof of (\ref{ApnK}).

{\bf The group $\mA_1^*$ and its commutants}.

\begin{theorem}\label{18Ma7}
Let $K$ be a field of characteristic zero. Then
\begin{enumerate}
\item $\mA_1^*= K^*\times (\CH \ltimes (1+F)^*)$, each unit $a$ of
$\mA_1$ is a unique product $a= \l \alpha (1+f)$ for some elements
$\l \in K^*$, $\alpha \in \CH$, and $f\in F$ such that $\det
(1+f)\neq 0$. \item $\mA_1^*= K^*\times (\CH \ltimes \GL_\infty
(K))$.
 \item The centre of the group $\mA_1^*$ is $K^*$.
 \item The commutant $\mA_1^{*(2)}:= [\mA_1^*, \mA_1^*]$ of the
 group $\mA_1^*$ is equal to $\SL_\infty (K):=\{ v\in (1+F)^*=
 M_\infty (K)\, | \, \det (v) = 1\}$, and $\mA_1/[\mA_1^*,
 \mA_1^*]\simeq K^*\times \CH \times K^*$.
 \item All the higher commutants $\mA_1^{*(i)}:= [\mA_1^*,
 \mA_1^{*(i-1)}]$, $i\geq 2$, are equal to $\mA_1^{*(2)}$.
\end{enumerate}
\end{theorem}

{\it Proof}. 1. By (\ref{ApnK}), $\mA_1'=K^*\times (\CH \ltimes
(1+F)^*)\subseteq \mA_1^*$. It suffices to show the reverse
inclusion. By (\ref{KHiA}), there is the exact sequence of groups
\begin{equation}\label{1FA1}
1\ra (1+F)^*\ra \mA_1^*\ra (\mA_1/F)^*\simeq \CA_1^*\ra 1
\end{equation}
which,  using again (\ref{KHiA}),  yields the inclusion
$\mA_1^*\subseteq \mA_1'$. The rest of statement 1 follows from
(\ref{1pFtn}).

 2. Statement 2 is equivalent to statement 1 since $(1+F)^* =
 \GL_\infty (K)$, see (\ref{1GLiK}).

3. Let $Z$ be the centre of the group $\mA_1^*$. Since
$K^*\subseteq Z$, we have
$$Z= Z\cap \mA_1^* = Z\cap (K^*\CH (1+F^*)= K^* (Z\cap \CH
(1+F)^*).$$ We have to show that $Z\cap \CH (1+F)^*= \{ 1\}$.
 Let $z= \alpha u\in Z\cap \CH (1+F)^*$ where $\alpha = \alpha
 (H)\in \CH$ and $u\in (1+F)^*$. It remains to show that $z=1$.
 $\beta z = z\beta$ for all $\beta \in \CH$ iff $\beta u = u\beta$
 for all $\beta \in \CH$ (since $\CH$ is an abelian group) iff
 $[\beta ] u = u [\beta ]$ (the equality of infinite matrices) for
 all $\beta = \beta (H)\in \CH$ where $[\beta ]$ is the diagonal
 matrix ${\rm diag} (\beta (1) , \beta (2) ,\ldots )$ iff $u$ is a {\em
 diagonal} matrix of $(1+F)^*$. The diagonal entries of the matrix $u$, say $u_i$,
  $i\in \N$, are elements of $K^*$ such that $u_i=1$ for all
  $i\geq d$ for some natural number $d=d(u)$.  For all distinct  $i,j\in
  \N$, $1+E_{ij}\in (1+F)^*$. Now, it follows from the equalities
  $$ z+\alpha (i+1) u_iE_{ij}= z(1+E_{ij})= (1+E_{ij})z= z+E_{ij}
  \alpha ( j+1) u_j$$
  that $\alpha (i+1) u_i= \alpha (j+1) u_j$ where $\alpha (i+1) :=
  \alpha (H)|_{H=i+1}$ (we have used that  $\alpha E_{ij}= \alpha (i+1)
  E_{ij}$ and
  $E_{ij}\alpha = E_{ij} \alpha (j+1)$). For all distinct natural
  numbers $i$ and $j$ such that $i,j>d$, we have $\alpha (i+1) =
  \alpha (j+1)$. This means that the function $\alpha (H)\in \CH$
  is a constant, i.e. $\alpha =1$. This gives $u_i= u_j$ for all
  $i,j\in \N$ such that $i\neq j$, i.e. all $u_i=1$. Therefore,
  $z=1$, as required.

4. The determinant can be extended from the subgroup $(1+F)^*$ of
$\mA_1^*$ to the whole group by the rule 
\begin{equation}\label{detA1*}
\det : \mA_1^*\ra K^*\times \CH \times K^*, \;\; \l \, \alpha \, u
\mapsto \l\, \alpha \, \det (u) := (\l, \alpha, \det (u)),
\end{equation}
where $\l \in K^*$, $\alpha \in \CH$, and $u\in (1+F)^*$. It turns
out that $\det $ is a group {\em epimorphism}. By the very
definition, $\det$ is a surjection. It remains to show that $\det
(aa') =\det (a) \det (a')$ for all $a:=\l\alpha u ,
a':=\l'\alpha'u'\in \mA_1^*$. This follows from the following
equality 
\begin{equation}\label{a1ua}
\det (\alpha^{-1} u \alpha ) = \det (u),\;\;  \alpha \in \CH, \;\;
u\in (1+F)^*.
\end{equation}
Indeed,
\begin{eqnarray*}
 \det (aa')&=& \det (\l\l'\alpha \alpha'\cdot (\alpha')^{-1} u\alpha'u')=  \l\l'\alpha \alpha'\cdot \det ((\alpha')^{-1} u\alpha')\det (u')\\
 &=&\l\l'\alpha \alpha'\cdot \det (u)\det (u')= \det (a) \det
 (a').
\end{eqnarray*}
The proof of (\ref{a1ua}):
\begin{eqnarray*}
\det  (\alpha^{-1} u \alpha ) &=& \det  ([\alpha^{-1}] u [\alpha ] )\;\;\;
({\rm where}\;\; [\alpha ] := {\rm diag}(\alpha (1), \alpha (2), \ldots )) \\
 &=&\det  ([\alpha^{-1}]_d u_d [\alpha ]_d )\;\;\; {\rm for\;\;
 all}\;\; d\gg 1\\
 &=& \det (u_d)= \det (u)  \;\;\; {\rm for\;\;
 all}\;\; d\gg 1
\end{eqnarray*}
where, for an infinite matrix $X=\sum_{i,j\in \N} x_{ij}E_{ij}$,
$X_d:= \sum_{0\leq i,j\leq d-1}x_{ij} E_{ij}$ is the sub-matrix of
$X$ of size $d\times d$.

The kernel of the epimorphism $\det$, (\ref{detA1*}), is
$\SL_\infty (K):= \{ u\in (1+F)^*= \GL_\infty (K)\, | \, \det (u)
=1\}$. In particular, $\SL_\infty (K)$ is a normal subgroup of
$\mA_1^*$ such that the factor group $\mA_1^* / \SL_\infty
(K)\simeq K^*\times \CH \times K^*$ is abelian. Hence, $[\mA_1^* ,
\mA_1^*]\subseteq \SL_\infty (K)$ and there is the short exact
sequence of groups 
\begin{equation}\label{sesA1}
1\ra \SL_\infty (K)\ra \mA_1^* \xrightarrow{\det } K^*\times \CH
\times K^*\ra 1.
\end{equation}
The group $\SL_\infty (K)$ is generated by the {\em transvections}
$t_{ij}(\l ):= 1+\l E_{ij}$, $\l \in K^*$, $i,j\in \N$, $i\neq j$.
Note that 
\begin{equation}\label{Htij}
[H, t_{ij}(\l )]= t_{ij}(\frac{i-j}{j+1}\, \l)
\end{equation}
where $[a,b]:=aba^{-1}b^{-1}$ is the commutator of two elements of
a group.  In more detail,
$$[H, t_{ij}(\l )]= Ht_{ij}(\l ) H^{-1}t_{ij}(\l)^{-1}= t_{ij}(\frac{i+1}{j+1}\,
\l) t_{ij}(-\l)= t_{ij}(\frac{i-j}{j+1}\, \l).
$$
Since $H\in \CH$, we have the inclusion $\SL_\infty (K) \subseteq
[\mA_1^*, \mA_1^*]$, i.e. $[\mA_1^*, \mA_1^*]=\SL_\infty (K)$.

5. Statement 5 follows from (\ref{Htij}) and the fact that the
group $\SL_\infty (K)$ is generated by the transvections.  $\Box $

{\bf An inversion formula for $u\in \mA_1^*$}. Let $K$ be a field
of characteristic zero. By Theorem \ref{18Ma7}.(1), each element
$u$ of $\mA_1^*$ can written as $u= \l a(1+f)$. The inverse
$(1+f)^{-1}$ can be found using Cramer's formula for the inverse
of matrix. Then $u^{-1} = \l^{-1} (1+f)^{-1} a^{-1}$. Let $f\in
K[x]$ be a given polynomial and $y\in K[x]$ is an unknown. Then
the integro-differential equation $uy=f$ can be solved explicitly:
$y= u^{-1}f$.

In contrast to differential operators on an affine line, in
general, the space of solutions for integro-differential operators
is infinite-dimensional: Example. $E_{ij} y=0$.

For an ideal $I$ of $\mA_n$ such that $I\neq \mA_n$, let $
(1+I)^*$ be the group of units of the multiplicative monoid $1+I$.

\begin{lemma}\label{I24Ma7}
Let $K$ be a commutative $\Q$-algebra, $I$ and $J$ be ideals of
$\mA_n$ which are distinct from $\mA_n$. Then
\begin{enumerate}
\item $\mA_n^*\cap (1+I)= (1+I)^*$.\item $(1+I)^*$ is a normal
subgroup of $\mA_n^*$.
\end{enumerate}
\end{lemma}

{\it Proof}. 1. The inclusion $\mA_n^*\cap (1+I)\supseteq (1+I)^*$
is obvious. To prove the reverse inclusion, let $1+a\in
\mA_n^*\cap (1+I)$ where $a\in I$, and let $(1+a)^{-1} = 1+b$ for
some $b\in \mA_n$. The equality $1= (1+a) (1+b)$ can be written as
$b= - a(1+b)\in I$, i.e. $1+a\in (1+I)^*$. This proves the reverse
inclusion.

2. For all $a\in \mA_n^*$, $a(1+I)a^{-1}= 1+aIa^{-1} = 1+I$, and
so $a(1+I)^*a^{-1} = a(\mA_n^* \cap (1+I))a^{-1} = a\mA_n^*a^{-1}
\cap a(1+I)a^{-1} = \mA_n^* \cap (1+I) = (1+I)^*$. Therefore,
$(1+I)^*$ is a normal subgroup of $\mA_n^*$.  $\Box $

Let $K$ be a field of characteristic zero. By (\ref{KHiA}), the
group homomorphism $\mA_n^*\ra ( \mA_n/ \ga_n)^*\simeq \CA_n^*$ is
an epimorphism. By Lemma \ref{I24Ma7}.(1), its kernel is
$\mA_n^*\cap (1+\ga_n) = (1+\ga_n)^*$, and we have the short exact
sequence of groups 
\begin{equation}\label{uFA1}
1\ra (1+\ga_n)^*\ra \mA_n^* \ra \CA_n^*\ra 1
\end{equation}
which together with (\ref{KHiA}) proves the first statement of the
next theorem.

\begin{theorem}\label{24Ma7}
Let $K$ be a field of characteristic zero. Then
\begin{enumerate}
\item $\mA_n^*= K^*\times (\CH_n\ltimes (1+\ga_n)^*)$. \item The
centre of the group $\mA_n^*$ is $K^*$.
\end{enumerate}
\end{theorem}

{\it Proof}. 2. Let $Z$ be the centre of the group $\mA_n^*$. Then
$K^*\subseteq Z$ and, by statement 1,
$$ Z= Z\cap \mA_n^* = Z\cap (K^*\CH_n(1+\ga_n)^*)= K^* (Z\cap
\CH_n (1+\ga_n)^*).$$ It remains to show that $Z':=Z\cap \CH_n
(1+\ga_n)^*=\{ 1\}$.

Let us show first that $Z'= Z\cap (1+\ga_n)^*$. Let $z= \v u\in
Z'$ for some $\v \in \CH_n$ and $u\in (1+\ga_n)^*$. It suffices to
show that $\v =1$. Note that, for each element $a\in \ga_n$, there
exists a natural number $c=c(a)$ such that $aE_{\alpha \beta} =
E_{\alpha \beta}a=0$  for all $\alpha, \beta \in \N^n$ such that
all $\alpha_i, \beta_i\geq c$. For short, we write $\alpha, \beta
\gg 0$. So, $uE_{\alpha, \beta } =E_{\alpha \beta}u= E_{\alpha
\beta}$ for all $\alpha , \beta \gg 0$. Note that $E_{\alpha
\beta}^2=0$ for all $\alpha \neq \beta$, and so $1+E_{\alpha
\beta}\in \mA_n^*$. Now, for all $\alpha, \beta \gg 0$ such that
$\alpha \neq \beta$, $z(1+E_{\alpha \beta}) = (1+E_{\alpha
\beta})z$ $ \Leftrightarrow$ $z+\v(\alpha +1) E_{\alpha \beta} =
z+E_{\alpha \beta} \v(\beta +1)$ $\Leftrightarrow$ $\v(\alpha +1)=
\v(\beta +1)$ $\Leftrightarrow$ $\v =1$, as required, where $\v
(\alpha +1)$ is the value of the function $\v = \v (H_1, \ldots ,
H_n)$ at $H_1=\alpha_1+1, \ldots , H_n= \alpha_n +1$. This proves
the claim.

So, it remains to show that $Z':= Z\cap ( 1+\ga_n)^*= \{ 1\}$. The
result is true for $n=1$ (Theorem \ref{18Ma7}.(3)).  So, we may
assume that $n\geq 2$. Consider the descending chain $\gf_1\supset
\cdots \supset \gf_i \supset \cdots \supset \gf_n\supset
\gf_{n+1}:=0$ of ideals
$$ \gf_i:= \sum_{1\leq j_1<\cdots <j_i\leq n } \gp_{j_1} \gp_{j_2}
\cdots \gp_{j_i}$$ of $\mA_n$. Note that $\gf_1= \ga_n$, $\gf_n =
F^{\t n }$, and
$$(\gp_{j_1}  \cdots \gp_{j_i}+\gf_{i+1})/ \gf_{i+1}
\simeq \gp_{j_1}  \cdots \gp_{j_i}/\gp_{j_1}  \cdots \gp_{j_i}\cap
\gf_{i+1}\simeq F^{\t i } \t \CA_{n-i} \simeq F\t \CA_{n-i} \simeq
M_\infty ( \CA_{n-i}).$$ In the proof of the above series of ring
isomorphisms without 1 we have used the facts that $F= M_\infty
(K)$ and $F^{\t i } \simeq F$. $$\gf_i/ \gf_{i+1} \simeq
\prod_{1\leq j_1<\cdots <j_i\leq n }(\gp_{j_1}  \cdots
\gp_{j_i}+\gf_{i+1})/ \gf_{i+1}\simeq (F^{\t i}\t \CA_{n-i})^{n
\choose i}\simeq M_\infty (\CA_{n-i})^{n\choose i}$$ are
isomorphisms of rings without 1. Consider the isomorphisms of
groups
\begin{eqnarray*}
 (1+\gf_i/ \gf_{i+1})^* &\simeq &(1+\prod_{1\leq j_1<\cdots <j_i\leq n
}(\gp_{j_1}  \cdots \gp_{j_i}+\gf_{i+1})/ \gf_{i+1})^*\\
&\simeq & \prod_{1\leq j_1<\cdots <j_i\leq n }(1+(\gp_{j_1}
\cdots \gp_{j_i}+\gf_{i+1})/ \gf_{i+1})^* \simeq {(1+F^{\t i}\t
\CA_{n-i})^*}^{n \choose i}\\
&\simeq & {(1+ M_\infty (\CA_{n-i}))^*}^{n\choose i}\simeq
\GL_\infty (\CA_{i-1})^{n\choose i}
\end{eqnarray*}
where $\GL_\infty (\CA_{i-1}):= (1+M_\infty (\CA_{n-i}))^*$. The
descending chain above yields the descending chain of normal
subgroups of $\mA_n^*$:
$$ (1+\gf_1)^* = (1+\ga_n)^*\supset \cdots \supset (1+\gf_i)^*\supset \cdots
\supset (1+\gf_n)^* = (1+F^{\t n })^*\supset (1+\gf_{n+1})^*= \{
1\}.$$ For each $i=1, \ldots , n$, there is the natural
homomorphism of groups $\v_i : (1+\gf_i)^*\ra
(1+\gf_i/\gf_{i+1})^*$, the kernel of which is $(1+\gf_{i+1})^*$.
So, $(1+\gf_i)^*/ (1+\gf_{i+1})^*$ is a subgroup of $(1+\gf_i/
\gf_{i+1})^*$. To prove that $Z'=\{ 1\}$ is equivalent to show
that $Z\cap (1+\gf_i)^* = \{ 1\}$ for all $i=1, \ldots , n$. To
prove this we use a {\em downward} induction at $i$ starting with
$i=n$. In this case, $\gf_n = F^{\t n }$,  and the fact that
$Z\cap (1+ F^{\t n })^*= \{ 1\}$ follows from the inclusion
$$ Z\cap (1+F^{\t n })^* \subseteq Z((1+ F^{\t n
})^*)=Z((1+M_\infty (K)^*)= \{ 1 \}$$ since $Z(M_\infty (K))= 0$.

Suppose that $i<n$ and $Z\cap (1+\gf_{i+1})^* = \{ 1\}$. Using a
downward induction on $i$ it remains to show that $Z_i:= Z\cap
(1+\gf_i)^* = \{ 1\}$. Note that in any ring elements $1+a$ and
$1+b$ commute iff the elements $a$ and $b$ commute. Using this
observation we see that,  for any ring $R$,  the centre
$Z(1+M_\infty (R))$ of the multiplicative monoid $1+M_\infty (R)$
is $1$ since $1$ is the only element of $1+M_\infty (K)$ that
commute with all the elements $1+E_{kl}$, $k\neq l$. All these
elements belong to the group $\GL_\infty (R):= (1+M_\infty
(R))^*$, and so it has the trivial centre 
\begin{equation}\label{ZGLIRt}
Z(\GL_\infty (R))=\{ 1\}.
\end{equation}
In particular, $Z(\GL_\infty (\CA_{n-i})^{n\choose i}) = \{ 1\}$.
For each subset $J= \{ j_1, \ldots , j_i\}$ of $\{ 1, \ldots ,
n\}$ that contains exactly $i$ elements, we have seen above that
$$ (\gp_{j_1}\cdots \gp_{j_i}+\gf_{i+1})/\gf_{i+1}\simeq F^{\t
i}\t \CA_{n-i}\simeq M_\infty (\CA_{n-i}).$$ For each matrix unit
$E_{\alpha \beta}=E_{\alpha_1\beta_1} (j_1) \cdots
E_{\alpha_i\beta_i}(j_i)\in \gp_{j_1} \cdots \gp_{j_i}$ where
$\alpha, \beta \in \N^i$ and $\alpha \neq \beta$, the elements
$1+E_{\alpha \beta}$ belongs to $(1+\gf_i)^*$ since $ E_{\alpha
\beta}^2=0$, and its image under the map $\v_i$ is equal to the
element $1+E_{\alpha \beta}+\gf_{i+1}$. This means that an element
$\v_i(z)$ commutes with all elements $1+E_{\alpha
\beta}+\gf_{i+1}$ for all possible choices of $J$, i.e. $\v_i (z)
\in Z(\GL_\infty (\CA_{n-i})^{n\choose i})=\{ 1\}$. This means
that $z\in Z\cap (1+\gf_{i+1})^*= \{ 1\}$, i.e. $z=1$, as
required. By induction, $\{ 1\} = Z\cap (1+\gf_1)^* = Z\cap
(1+\ga_n)^*$. This proves that $Z(\mA_n^*)= K^*$. $\Box $

{\bf The subgroup $(1+\ga_n)^*$ of $\mA_n^*$}. Let $K$ be a
commutative $\Q$-algebra. For each $i=1, \ldots , n$,
\begin{eqnarray*}
(\gp_i +\gf_2)/ \gf_2&\simeq & \gp_i/ \gp_i\cap \gf_2\simeq
\mA_{i-1} \t F \t \mA_{n-i} / ( \ga_{i-1} \t F \t \mA_{n-i}+
\mA_{i-1} \t F \t \ga_{n-i})\\
&\simeq & \mA_{i-1} / \ga_{i-1} \t F \t \mA_{n-i} /
\ga_{n-i}\simeq \CA_{i-1}\t F \t \CA_{n-i} \simeq M_\infty (
\CA_{n-1})
\end{eqnarray*}
is the series of isomorphisms of rings without 1. The factor ring
without 1 $$\ga_n/ \gf_2\simeq \prod_{i=1}^n (\gp_i + \gf_2) /
\gf_2\simeq \prod_{i=1}^n \CA_{i-1} \t F \t \CA_{n-i}\simeq
M_\infty (\CA_{n-1})^n$$ is the direct product of its subrings
without 1. It is a semi-simple $\mA_n$-bimodule where $(\gp_i
+\gf_2) / \gf_2$, $1\leq i \leq n$, are the {\em simple} isotypic
components of the $\mA_n$-bimodule $\ga_n/ \gf_2$. Fix a section
$s: \ga_n/ \gf_2\ra \ga_n$ of the $K$-module epimorphism $\ga_n\ra
\ga_n /\gf_2$, $a\mapsto a+\gf_2$. Then $\ga_n = \im (s) \oplus
\gf_2$ is a direct sum  of $K$-submodules. Using the $K$-basis
$B_n$  for $\mA_n$ considered in the proof of Corollary
\ref{15Ma7} one can easily find such a section which even
satisfies the additional property that $\im (s) $ is a free
$K$-module.

Consider the ring $ K\oplus (\ga_n / \gf_2)$ with 1 and the
subgroup $(1+ (\ga_n /\gf_2))^*$ of its group $(K\oplus (\ga_n /
\gf_2))^*$ of units. There are canonical group isomorphisms
\begin{eqnarray*}
 (1+ \ga_n /\gf_2)^* & \simeq & (1+\prod_{i=1}^n (\gp_i +\gf_2)/\gf_2)^* \\
 &\simeq & \prod_{i=1}^n ( 1+(\gp_i +\gf_2) / \gf_2)^*, \;\; 1+\sum_{i=1}^n p_i\mapsto \prod_{i=1}^n (1+p_i), \\
 &\simeq & \prod_{i=1}^n (1+M_\infty (\CA_{n-1}))^* =
 \prod_{i=1}^n \GL_\infty ( \CA_{n-1})= \GL_\infty (\CA_{n-1})^n.
\end{eqnarray*}
We have the group monomorphism
$$ (1+\ga_n)^* / (1+\gf_2)^* \ra (1+ \ga_n / \gf_2)^*, \;\;
(1+a)(1+\gf_2)^*\mapsto 1+a+\gf_2.$$

{\bf An invertibility criterion}.  The next theorem is a criterion
of when an element of the monoid $1+\ga_n$ is invertible.
\begin{theorem}\label{23Ma7}
Let $K$ be a commutative $\Q$-algebra, $a\in \ga_n$. Then $1+a\in
(1+\ga_n)^*$ iff
\begin{enumerate}
\item $1+a+\gf_2\in (1+\ga_n/\gf_2)^*\, (\simeq \GL_\infty
(\CA_{n-1})^n)$, and \item $a+(1+a)c\in (1+a) \gf_2$ and
$a+c(1+a)\in \gf_2 (1+a)$ where $c:= s((1+a+\gf_2)^{-1} -1)$ (the
value of the section $s:\ga_n/\gf_2\ra \ga_n$ at the element
$(1+a+\gf_2)^{-1} -1\in \ga_n/\gf_2$).
\end{enumerate}
Suppose that conditions 1 and 2 hold and $a+(1+a)c=(1+a)r$ (resp.
$a+c(1+a)=l(1+a)$) for some $r\in \gf_2$ (resp. $l\in \gf_2$) then
$a^{-1}= 1+c-r$ (resp. $a^{-1} = 1+c-l$).
\end{theorem}

{\it Proof}. $(\Rightarrow )$ Suppose that $(1+a)\in (1+\ga_n)^*$.
Due to the group homomorphism $(1+\ga_n)^*\ra (1+\ga_n/\gf_2)^*$,
we have $1+a\mapsto 1+a+\gf_2\in (1+\ga_n/\gf_2)^*$, i.e. the
first condition holds. $(1+\ga_n)^*\ni (1+a)^{-1} = 1+b$ for some
element $b\in \ga_n= \im (s) \oplus \gf_2$ which can be written as
$b=c+d$ where $c:= s((1+a+\gf_2)^{-1} -1)\in \im (s)$ and $d:=
b-c\in \gf_2$. The equalities $(1+a) (1+c+d) =1 $ and $(1+c+d)
(1+a)=1$ can be rewritten as follows $a+(1+a)c= -(1+a) d\in (1+a)
\gf_2$ and $a+c(1+a) = - d(1+a) \in \gf_2 (1+a)$, and so the
second statement holds.

$( \Leftarrow )$ Suppose that conditions 1 and 2 hold, we have to
show that $1+a\in (1+\ga_n)^*$, i.e. the element $1+a$ has a left
and a right inverse. Condition 2 can be     written  as follows
$a+(1+a) c = (1+a) r$ and $a+c(1+a) = l (1+a)$ for some elements
$r,l\in \gf_2$. These two equalities can we rewritten as $(1+a)
(1+c-r)=1$ and $ (1+c-l)(1+a) = 1$. This means that $1+a\in
(1+\ga_n)^*$ and $a^{-1} =  1+c-r= 1+c-l$.  $\Box $

{\bf An inversion formula for $u\in K^*\times (\CH_n\ltimes
(1+F^{\t n })^*)$}. Let $K$ be a field of characteristic zero. By
Theorem \ref{24Ma7}.(1), each element $u\in \mA_n^*$ is a unique
product $u = \l h (1+a)$ for some $\l \in K^*$, $h\in \CH_n$, and
$a\in \ga_n$ such that $1+a\in (1+\ga_n)^*$. Clearly, $u^{-1}=
\l^{-1} (1+a)^{-1} h^{-1}$. So, to write down explicitly an
inversion formula for $u$ boils down to finding $(1+a)^{-1}$. As a
first step, one should know an inversion formula for elements of
$\GL_\infty (\CA_{n-1})$ which is not obvious how to do, at the
moment. It should not be entirely trivial since as a result one
would have a formula for solutions of all invertible
integro-differential equations (for all $n\geq 2$). Nevertheless,
for elements $u= \l h (1+a)\in K^*\times (\CH_n\ltimes (1+F^{\t n
})^*)$ one {\em can} write down the inversion formula exactly in
the same manner as in the case $n=1$. Hence, one obtains
explicitly solutions to the equation $uy=f$ where $f\in P_n$.
Elements of the group $K^*\times (\CH_n\ltimes (1+F^{\t n })^*)$
are called {\em minimal integro-differential operators}.

Department of Pure Mathematics

University of Sheffield

Hicks Building

Sheffield S3 7RH

UK

email: v.bavula@sheffield.ac.uk

\end{document}